\documentclass[thmsa]{article}
\usepackage{amsfonts}
\usepackage{amsmath,amsthm,amssymb}
\usepackage{graphicx}
\usepackage[noadjust]{cite}

\newtheorem{theorem}{Theorem}
\newtheorem{corollary}[theorem]{Corollary}
\newtheorem{lemma}[theorem]{Lemma}
\newtheorem{proposition}[theorem]{Proposition}

\newtheorem{remark}[theorem]{Remark}

\def\proof{\noindent{{\bf Proof. }}}
\def\endproof{ $\blacksquare$}

\DeclareMathOperator*{\esssup}{ess\,sup}
\DeclareMathOperator*{\essinf}{ess\,inf}
\DeclareMathOperator*{\linspan}{span}

\makeatletter
\let\@fnsymbol\@alph
\makeatother

\begin{document}

\title{Compactness and existence results in weighted Sobolev spaces of radial
functions. Part II: Existence}

\author{Marino Badiale\thanks{Dipartimento di Matematica ``Giuseppe Peano'', Universit\`{a} degli Studi di
Torino, Via Carlo Alberto 10, 10123 Torino, Italy. 
e-mails: \texttt{marino.badiale@unito.it}, \texttt{michela.guida@unito.it}}
\textsuperscript{,}\thanks{Partially supported by the PRIN2012 grant ``Aspetti variazionali e
perturbativi nei problemi di.renziali nonlineari''.}
\ -\ Michela Guida\textsuperscript{a,}\thanks{Member of the Gruppo Nazionale di Alta Matematica (INdAM).}
\ -\ Sergio Rolando\thanks{Dipartimento di Matematica e Applicazioni, Universit\`{a} di Milano-Bicocca,
Via Roberto Cozzi 53, 20125 Milano, Italy. e-mail: \texttt{sergio.rolando@unito.it}}%
\ \textsuperscript{,c}
}
\date{}
\maketitle

\begin{abstract}
We apply the compactness results obtained in \cite{BGR-p1} to prove
existence and multiplicity results for finite energy solutions to the
nonlinear elliptic equation 
\[
-\triangle u+V\left( \left| x\right| \right) u=g\left( \left| x\right|
,u\right) \quad \textrm{in }\Omega \subseteq \mathbb{R}^{N},\ N\geq 3, 
\]
where $\Omega $ is a radial domain (bounded or unbounded) and $u$ satisfies 
$u=0$ on $\partial \Omega $ if $\Omega \neq \mathbb{R}^{N}$ and $u\rightarrow 0$
as $\left| x\right| \rightarrow \infty $ if $\Omega $ is unbounded. The
potential $V$ may be vanishing or unbounded at zero or at infinity and the
nonlinearity $g$ may be superlinear or sublinear. If $g$ is sublinear, the
case with $g\left( \left| \cdot \right| ,0\right) \neq 0$ is also
considered.\bigskip

\smallskip

\medskip
\end{abstract}

\section{Introduction and main results}

In this paper we study the existence and multiplicity of radial solutions to the following problem: 
\begin{equation}
\left\{ 
\begin{array}{l}
-\triangle u+V\left( \left| x\right| \right) u=g\left( \left| x\right|,u\right) 
\quad \text{in }\Omega \medskip \\ 
u\in D_{0}^{1,2}(\Omega ),\quad \int_{\Omega }V\left( \left| x\right|\right) u^{2}dx<\infty
\end{array}
\right.  \tag*{$\left( P\right) $}
\end{equation}
where $\Omega \subseteq \mathbb{R}^{N}$, $N\geq 3$, is a spherically symmetric
domain (bounded or unbounded), $D_{0}^{1,2}(\Omega )$ is the usual Sobolev
space given by the completion of $C_{\mathrm{c}}^{\infty }(\Omega )$ with
respect to the $L^{2}$ norm of the gradient and the potential $V$ satisfies
the following basic assumption, where 
$\Omega _{\mathrm{r}}:=\left\{ \left|x\right| >0:x\in \Omega \right\}$:

\begin{itemize}
\item[$\left( \mathbf{V}\right) $]  
$V:\Omega _{\mathrm{r}}\rightarrow\left[ 0,+\infty \right) $ is a measurable function such that 
$V\in L^{1}\left( r_{1},r_{2}\right) $ for some interval $\left(r_{1},r_{2}\right) \subseteq \Omega _{\mathrm{r}}.$
\medskip
\end{itemize}

\noindent More precisely, we define the space 
\begin{equation}
H_{0,V}^{1}\left( \Omega \right) :=\left\{ u\in D_{0}^{1,2}\left( \Omega
\right) :\int_{\Omega }V\left( \left| x\right| \right) u^{2}dx<\infty
\right\}  \label{H^1_0,V :=}
\end{equation}
(which is nonzero by assumption $\left( \mathbf{V}\right) $) and look for
solutions in the following weak sense: we name \textit{solution} to problem $%
\left( P\right) $ any $u\in H_{0,V}^{1}\left( \Omega \right) $ such that 
\begin{equation}
\int_{\Omega }\nabla u\cdot \nabla h\,dx+\int_{\Omega }V\left( \left|
x\right| \right) uh\,dx=\int_{\Omega }g\left( \left| x\right| ,u\right)
h\,dx\qquad \text{for all }h\in H_{0,V}^{1}\left( \Omega \right) .
\label{weak solution}
\end{equation}
Of course, we will say that a solution is \textit{radial} if it is invariant
under the action on $H_{0,V}^{1}\left( \Omega \right) $ of the orthogonal
group of $\mathbb{R}^{N}$. 

By well known arguments, problem $\left( P\right) $ is a model for the
stationary states of reaction diffusion equations in population dynamics
(see e.g. \cite{Fife}). Moreover, its nonnegative weak solutions lead to
special solutions (\textit{solitary waves} and \textit{solitons}) for
several nonlinear field theories, such as nonlinear Schr\"{o}dinger (or
Gross-Pitaevskii) and Klein-Gordon equations, which arise in many branches
of mathematical physics, such as nonlinear optics, plasma physics, condensed
matter physics and cosmology (see e.g. \cite{BFmonograph,YangY}). In this
respect, since the early studies of \cite
{Beres-Lions,Floer-Wein,Rabi92,Strauss}, problem $\left( P\right) $ has been
massively addressed in the mathematical literature, recently focusing on the
case with $\Omega =\mathbb{R}^{N}$ and $V$ possibly vanishing at infinity, that
is, $\liminf_{\left| x\right| \rightarrow \infty }V\left( \left| x\right|
\right) =0$ (some first results on such a case can be found in \cite
{Ambr-Fel-Malch,BR,Be-Gr-Mic,Be-Gr-Mic-2}; for more recent bibliography, see e.g. \cite{Alves-Souto-13,BGRnonex,BR-TMA,Catrina-nonex,SuTian12,Chen13,BonMerc11,Yang-Li-14,Su12,Zhang13,Deng-Peng-Pi-14,Su-Wang-Will-p}
and the references therein).

Here we study problem $\left( P\right) $ under assumptions that, together
with $\left( \mathbf{V}\right) $, allow $V\left( r\right) $ to be singular
at some points (including the origin if $\Omega $ is a ball), or vanishing
as $r\rightarrow +\infty $ (if $\Omega $ is unbounded), or both. Also the
case of $V=0$, or $V$ compactly supported, or $V$ vanishing in a neighbourhood of the origin, 
will be encompassed by our results.
As concerns the nonlinearity, we will mainly focus
on the following model case (see Section \ref{SEC: gen-results} for more
general results): 
\begin{equation}
g\left( \left| x\right| ,u\right) =K\left( \left| x\right| \right) f\left(
u\right)  \label{P0}
\end{equation}
where $f$ and the potential $K$ satisfy the following basic assumptions:

\begin{itemize}
\item[$\left( \mathbf{K}\right) $]  $K:\Omega _{\mathrm{r}}\rightarrow
\left( 0,+\infty \right) $ is a measurable function such that $K\in L_{%
\mathrm{loc}}^{s}\left( \Omega _{\mathrm{r}}\right) $ for some $s>\frac{2N}{%
N+2};$

\item[$\left( \mathbf{f}\right) $]  $f:\mathbb{R}\rightarrow \mathbb{R}$ is
continuous and such that $f\left( 0\right) =0$.
\end{itemize}

\noindent Both the cases of $f$ superlinear and sublinear will be studied.
For sublinear $f$, we will also deal with an additional forcing term, i.e.,
with nonlinearities of the form: 
\begin{equation}
g\left( \left| x\right| ,u\right) =K\left( \left| x\right| \right) f\left(
u\right) +Q\left( \left| x\right| \right) .  \label{PQ}
\end{equation}
Problem $\left( P\right) $ with such $g$'s will be denoted by $\left(
P_{Q}\right) $, so that, accordingly, $\left( P_{0}\right) $ will indicate
problem $\left( P\right) $ with $g$ given by (\ref{P0}).
\smallskip


Besides hypotheses $\left( \mathbf{V}\right) $, $\left( \mathbf{K}\right) $
and $\left( \mathbf{f}\right) $, which will be always tacitly assumed in
this section, the potentials $V$ and $K$ will satisfy suitable combinations
of the following conditions:

\begin{itemize}
\item[$\left( \mathbf{VK}_{0}\right) $]  $\exists \alpha _{0}\in \mathbb{R}$
and $\exists \beta _{0}\in \left[ 0,1\right] $ such that 
\[
\esssup_{r\in \left( 0,R_{0}\right) }\frac{K\left( r\right) }{%
r^{\alpha _{0}}V\left( r\right) ^{\beta _{0}}}<+\infty \quad \text{for some }%
R_{0}>0;
\]

\item[$\left( \mathbf{VK}_{\infty }\right) $]  $\exists \alpha _{\infty }\in 
\mathbb{R}$ and $\exists \beta _{\infty }\in \left[ 0,1\right] $ such that 
\[
\esssup_{r>R_{\infty }}\frac{K\left( r\right) }{r^{\alpha _{\infty
}}V\left( r\right) ^{\beta _{\infty }}}<+\infty \quad \text{for some }%
R_{\infty }>0;
\]

\item[$\left( \mathbf{V}_{0}\right) $]  $\exists \gamma _{0}>2$ such that $%
\essinf\limits_{r\in \left( 0,R_{0}\right) }r^{\gamma _{0}}V\left(
r\right) >0$ for some $R_{0}>0;$

\item[$\left( \mathbf{V}_{\infty }\right) $]  $\exists \gamma _{\infty }<2$
such that $\essinf\limits_{r>R_{\infty }}r^{\gamma _{\infty
}}V\left( r\right) >0$ for some $R_{\infty }>0.$
\end{itemize}

\noindent We mean that $V\left( r\right) ^{0}=1$ for every $r$, so that conditions 
$\left( \mathbf{VK}_{0}\right) $ and $\left( \mathbf{VK}_{\infty}\right) $ 
will also make sense if $V\left( r\right) =0$ for $r<R_{0}$ or $r>R_{\infty }$, 
with $\beta _{0}=0$ or $\beta _{\infty }=0$ respectively.
\smallskip

Concerning the nonlinearity, our existence results rely on suitable
combinations of the following assumptions:

\begin{itemize}
\item[$\left( \mathbf{f}_{1}\right) $]  $\exists q_{1},q_{2}>1$ such that 
\[
\sup_{t>0}\,\frac{\left| f\left( t\right) \right| }{\min \left\{
t^{q_{1}-1},t^{q_{2}-1}\right\} }<+\infty ;
\]

\item[$\left( \mathbf{F}_{1}\right) $]  $\exists \theta >2$ and $\exists
t_{0}>0$ such that $0\leq \theta F\left( t\right) \leq f\left( t\right) t$
for all $t\geq 0$ and $F\left( t_{0}\right) >0;$

\item[$\left( \mathbf{F}_{2}\right) $]  $\exists \theta >2$ and $\exists
t_{0}>0$ such that $0<\theta F\left( t\right) \leq f\left( t\right) t$ for
all $t\geq t_{0};$

\item[$\left( \mathbf{F}_{3}\right) $]  $\exists \theta <2$ such that $%
\displaystyle%
\liminf_{t\rightarrow 0^{+}}\frac{F\left( t\right) }{t^{\theta }}>0.$
\end{itemize}

\noindent Here and in the following, we denote $F\left( t\right)
:=\int_{0}^{t}f\left( s\right) ds$. Observe that the double-power growth
condition $\left( \mathbf{f}_{1}\right) $ with $q_{1}\neq q_{2}$ is more
stringent than the following and more usual single-power one:

\begin{itemize}
\item[$\left( \mathbf{f}_{2}\right) $]  $\exists q>1$ such that 
\[
\sup_{t>0}\,\frac{\left| f\left( t\right) \right| }{t^{q-1}}<+\infty 
\]
\end{itemize}

\noindent (the former implies the latter for $q=q_{1}$, $q=q_{2}$ and every $%
q$ in between). On the other hand, $\left( \mathbf{f}_{1}\right) $ does not
require $q_{1}\neq q_{2}$, so that it is actually equivalent to $\left( 
\mathbf{f}_{2}\right) $ as long as one can take $q_{1}=q_{2}$ (cf. Remarks 
\ref{RMK:RNsuper}.\ref{RMK:RNsuper: single} and \ref{RMK:RNsub}.\ref
{RMK:RNsub: single}).\medskip 

In order to state our existence results for superlinear nonlinearities, we
introduce the following notation. For $\alpha ,\gamma \in \mathbb{R}$ and $%
\beta \in \left[ 0,1\right] $, we define the function 
\begin{equation}
q_{**}\left( \alpha ,\beta ,\gamma \right) :=2\frac{2\alpha +\left( 1-2\beta
\right) \gamma +2\left( N-1\right) }{2\left( N-1\right) -\gamma }\qquad 
\text{if }\gamma \neq 2N-2.  \label{q** :=}
\end{equation}
Then, for $\gamma \geq 2$, we set 
\[
\underline{\alpha }\left( \beta ,\gamma \right) :=\left\{ 
\begin{array}{ll}
-\left( 1-\beta \right) \gamma ~~\smallskip & \text{if }2\leq \gamma \leq
2N-2 \\ 
-\infty & \text{if }\gamma >2N-2,
\end{array}
\right. 
\]
\[
\underline{q}\left( \alpha ,\beta ,\gamma \right) :=\left\{ 
\begin{array}{ll}
2\smallskip & \text{if }2\leq \gamma \leq 2N-2\text{~and }\alpha >\underline{%
\alpha }\left( \beta ,\gamma \right) \\ 
\max \left\{ 2,q_{**}\left( \alpha ,\beta ,\gamma \right) \right\} ~~ & 
\text{if }\gamma >2N-2\text{~and }\alpha >\underline{\alpha }\left( \beta
,\gamma \right) ,
\end{array}
\right. 
\]
and 
\[
\overline{q}\left( \alpha ,\beta ,\gamma \right) :=\left\{ 
\begin{array}{ll}
q_{**}\left( \alpha ,\beta ,\gamma \right) ~~\smallskip & \text{if }2\leq
\gamma <2N-2\text{~and }\alpha >\underline{\alpha }\left( \beta ,\gamma
\right) \\ 
+\infty & \text{if }\gamma \geq 2N-2\text{~and }\alpha >\underline{\alpha }%
\left( \beta ,\gamma \right) .
\end{array}
\right. 
\]

\begin{theorem}
\label{THM:RNsuper}Let $\Omega =\mathbb{R}^{N}$. Assume that $f$ satisfies $%
\left( \mathbf{F}_{1}\right) $, or that $K\left( \left| \cdot \right|
\right) \in L^{1}(\mathbb{R}^{N})$ and $f$ satisfies $\left( \mathbf{F}%
_{2}\right) $. Assume furthermore that $\left( \mathbf{VK}_{0}\right) $, $%
\left( \mathbf{VK}_{\infty }\right) $ and $\left( \mathbf{f}_{1}\right) $
hold with 
\begin{equation}
\alpha _{0}>\underline{\alpha },\qquad \underline{q}<q_{1}<\overline{q}%
,\qquad q_{2}>\max \left\{ 2,q_{**}\right\} ,  \label{THM:RNsuper: ineqs}
\end{equation}
where 
\[
\underline{\alpha }=\underline{\alpha }\left( \beta _{0},2\right) ,\quad 
\underline{q}=\underline{q}\left( \alpha _{0},\beta _{0},2\right) ,\quad 
\overline{q}=\overline{q}\left( \alpha _{0},\beta _{0},2\right) \quad \text{%
and}\quad q_{**}=q_{**}\left( \alpha _{\infty },\beta _{\infty },2\right) .
\]
Then problem $\left( P_{0}\right) $ has a nonnegative radial solution $u\neq
0$. If $V$ also satisfies $\left( \mathbf{V}_{0}\right) $, then we can take $%
\underline{\alpha }=\underline{\alpha }\left( \beta _{0},\gamma _{0}\right) $%
, $\underline{q}=\underline{q}\left( \alpha _{0},\beta _{0},\gamma
_{0}\right) $ and $\overline{q}=\overline{q}\left( \alpha _{0},\beta
_{0},\gamma _{0}\right) $. If $V$ also satisfies $\left( \mathbf{V}_{\infty
}\right) $, then we can take $q_{**}=q_{**}\left( \alpha _{\infty },\beta
_{\infty },\gamma _{\infty }\right) $.
\end{theorem}

\begin{remark}
\label{RMK:RNsuper}\quad 

\begin{enumerate}
\item  The inequality $\underline{q}<\overline{q}$ is not an assumption in (%
\ref{THM:RNsuper: ineqs}) (even in the cases with assumptions $\left( 
\mathbf{V}_{0}\right) $, $\left( \mathbf{V}_{\infty }\right) $), since it is
ensured by the condition $\alpha _{0}>\underline{\alpha }$.

\item  \label{RMK:RNsuper: improve}For $\beta \in \left[ 0,1\right] $ fixed, 
$\underline{\alpha }\left( \beta ,\gamma \right) $ is left-continuous and
decreasing in $\gamma \geq 2$ (as a real extended function) and one can
check that $\underline{q}\left( \alpha ,\beta ,\gamma \right) $ and $%
\overline{q}\left( \alpha ,\beta ,\gamma \right) $, defined on the set $%
\left\{ \left( \alpha ,\gamma \right) :\gamma \geq 2,\,\alpha >\underline{%
\alpha }\left( \beta ,\gamma \right) \right\} $, are continuous and
respectively decreasing and increasing, both in $\gamma $ for $\alpha $
fixed and in $\alpha $ for $\gamma $ fixed ($\overline{q}$ is continuous and
increasing as a real extended valued function). Similarly, $\max \left\{
2,q_{**}\left( \alpha ,\beta ,\gamma \right) \right\} $ is increasing and
continuous both in $\gamma \leq 2$ for $\alpha \in \mathbb{R}$ fixed and in $%
\alpha \in \mathbb{R}$ for $\gamma \leq 2$ fixed.

Therefore, thanks to such monotonicities in $\gamma $, Theorem \ref
{THM:RNsuper} actually improves under assumption $\left( \mathbf{V}%
_{0}\right) $, or $\left( \mathbf{V}_{\infty }\right) $, or both.

Moreover, by both monotonicity and continuity (or left-continuity) in $%
\alpha $ and $\gamma $, the theorem is also true if we replace $\alpha
_{0},\alpha _{\infty },\gamma _{0},\gamma _{\infty }$ in $\underline{\alpha }%
,\underline{q},\overline{q},q_{**}$ with $\overline{\alpha }_{0},\underline{%
\alpha }_{\infty },\overline{\gamma }_{0},\underline{\gamma }_{\infty }$,
where 
\[
\overline{\alpha }_{0}:=\sup \left\{ \alpha _{0}:\esssup_{r\in
\left( 0,R_{0}\right) }\frac{K\left( r\right) }{r^{\alpha _{0}}V\left(
r\right) ^{\beta _{0}}}<+\infty \right\} ,\quad \underline{\alpha }_{\infty
}:=\inf \left\{ \alpha _{\infty }:\esssup_{r>R_{\infty }}\frac{%
K\left( r\right) }{r^{\alpha _{\infty }}V\left( r\right) ^{\beta _{\infty }}}%
<+\infty \right\} ,
\]
\[
\overline{\gamma }_{0}:=\sup \left\{ \gamma _{0}:\essinf%
\limits_{r\in \left( 0,R_{0}\right) }r^{\gamma _{0}}V\left( r\right)
>0\right\} ,\quad \underline{\gamma }_{\infty }:=\inf \left\{ \gamma
_{\infty }:\essinf\limits_{r>R_{\infty }}r^{\gamma _{\infty
}}V\left( r\right) >0\right\} .
\]
This is consistent with the fact that $\left( \mathbf{VK}_{0}\right) ,\left( 
\mathbf{V}_{0}\right) $ and $\left( \mathbf{VK}_{\infty }\right) ,\left( 
\mathbf{V}_{\infty }\right) $ still hold true if we respectively lower $%
\alpha _{0},\gamma _{0}$ and raise $\alpha _{\infty },\gamma _{\infty }$.

\item  \label{RMK:RNsuper: single}Theorem \ref{THM:RNsuper} also concerns
the case of power-like nonlinearities, since the exponents $q_{1}$ and $q_{2}
$ need not to be different in $\left( \mathbf{f}_{1}\right) $ and one can
take $q_{1}=q_{2}$ as soon as $\max \left\{ 2,q_{**}\right\} <\overline{q}$.
For example, this is always the case when $\left( \mathbf{V}_{0}\right) $
holds with $\gamma _{0}\geq 2N-2$ (which gives $\overline{q}=+\infty $), or
when $\alpha _{\infty }\leq 2\left( \beta _{\infty }-1\right) $ (which
implies $\max \left\{ 2,q_{**}\left( \alpha _{\infty },\beta _{\infty
},2\right) \right\} =2$).
\end{enumerate}
\end{remark}

The Dirichlet problem in bounded ball domains or exterior spherically
symmetric domains can be reduced to the problem in $\Omega =\mathbb{R}^{N}$ by
suitably modifying the potentials $V$ and $K$ (see Section \ref{SEC: pfs} below).
Hence, by the same arguments yielding Theorem \ref{THM:RNsuper}, we will
also get the following results.

\begin{theorem}
\label{THM:Bsuper}Let $\Omega $ be a bounded ball. Assume that $f$ satisfies 
$\left( \mathbf{F}_{1}\right) $, or that $K\left( \left| \cdot \right|
\right) \in L^{1}(\Omega )$ and $f$ satisfies $\left( \mathbf{F}_{2}\right) $%
. Assume furthermore that $\left( \mathbf{VK}_{0}\right) $ and $\left( 
\mathbf{f}_{2}\right) $ hold with 
\begin{equation}
\alpha _{0}>\underline{\alpha }\quad \text{and}\quad \underline{q}<q<%
\overline{q},\qquad \text{where\quad }\underline{\alpha }=\underline{\alpha }%
\left( \beta _{0},2\right) ,~\underline{q}=\underline{q}\left( \alpha
_{0},\beta _{0},2\right) ,~\overline{q}=\overline{q}\left( \alpha _{0},\beta
_{0},2\right) .  \label{THM:Bsuper: ineqs}
\end{equation}
Then problem $\left( P_{0}\right) $ has a nonnegative radial solution $u\neq
0$. If $V$ also satisfies $\left( \mathbf{V}_{0}\right) $, then we can take $%
\underline{\alpha }=\underline{\alpha }\left( \beta _{0},\gamma _{0}\right) $%
, $\underline{q}=\underline{q}\left( \alpha _{0},\beta _{0},\gamma
_{0}\right) $ and $\overline{q}=\overline{q}\left( \alpha _{0},\beta
_{0},\gamma _{0}\right) $.
\end{theorem}

\begin{theorem}
\label{THM:BCsuper}Let $\Omega $ be an exterior radial domain. Assume that $f
$ satisfies $\left( \mathbf{F}_{1}\right) $, or that $K\left( \left| \cdot
\right| \right) \in L^{1}(\Omega )$ and $f$ satisfies $\left( \mathbf{F}%
_{2}\right) $. Assume furthermore that $\left( \mathbf{VK}_{\infty }\right) $
and $\left( \mathbf{f}_{2}\right) $ hold with 
\[
q>\max \left\{ 2,q_{**}\right\} ,\qquad \text{where\quad }%
q_{**}=q_{**}\left( \alpha _{\infty },\beta _{\infty },2\right) .
\]
Then problem $\left( P_{0}\right) $ has a nonnegative radial solution $u\neq
0$. If $V$ also satisfies $\left( \mathbf{V}_{\infty }\right) $, then we can
take $q_{**}=q_{**}\left( \alpha _{\infty },\beta _{\infty },\gamma _{\infty
}\right) $.
\end{theorem}

For dealing with the sublinear case, we need some more notation. For $\alpha
,\gamma \in \mathbb{R}$ and $\beta \in \left[ 0,1\right) $, we define the
following functions: 
\begin{equation}
\alpha _{1}\left( \beta ,\gamma \right) :=-\left( 1-\beta \right) \gamma
,\qquad \alpha _{2}\left( \beta \right) :=-\left( 1-\beta \right) N,\qquad
\alpha _{3}\left( \beta ,\gamma \right) :=-\frac{\left( 1-2\beta \right)
\gamma +N}{2}  \label{alphas :=}
\end{equation}
and 
\begin{equation}
q_{*}\left( \alpha ,\beta ,\gamma \right) :=2\frac{\alpha -\gamma \beta +N}{%
N-\gamma }\qquad \text{if }\gamma \neq N.  \label{q* :=}
\end{equation}
Then, for $\gamma \geq 2$, we set 
\[
q_{0}\left( \alpha ,\beta ,\gamma \right) :=\left\{ 
\begin{array}{ll}
\max \left\{ 1,2\beta \right\} \smallskip & \text{if }2\leq \gamma \leq N%
\text{~and }\alpha \geq \alpha _{1}\left( \beta ,\gamma \right) \\ 
\max \left\{ 1,2\beta ,q_{*}\left( \alpha ,\beta ,\gamma \right) \right\} ~~
& \text{if }\gamma >N\text{~and }\alpha \geq \alpha _{1}\left( \beta ,\gamma
\right) .
\end{array}
\right. 
\]
In contrast with the superlinear case, we divide our existence results into
two theorems, essentially according as assumption $\left( \mathbf{VK}%
_{0}\right) $ holds with $\alpha _{0}$ large enough with respect to $\beta
_{0}$ (and $\gamma _{0}$, if $\left( \mathbf{V}_{0}\right) $ holds), or not
(cf. Remark \ref{RMK:RNsub}.\ref{RMK:RNsub: useless}): in the first case, we
need only require that $f$ grows as a single power; in the second case, we
assume the double-power growth condition $\left( \mathbf{f}_{1}\right) $,
which, however, may still reduce to a single-power one in particular cases
of exponents $\alpha _{0},\beta _{0},\gamma _{0},\alpha _{\infty },\beta
_{\infty },\gamma _{\infty }$ (cf. Remark \ref{RMK:RNsub}.\ref{RMK:RNsub:
single}).

\begin{theorem}
\label{THM:RNsub1}Let $\Omega =\mathbb{R}^{N}$ and let $Q\in L^{2}(\mathbb{R}%
_{+},r^{N+1}dr)$, $Q\geq 0$. Assume that $f$ satisfies $\left( \mathbf{F}%
_{3}\right) $, or that $Q$ does not vanish almost everywhere in $\left(
r_{1},r_{2}\right) $. Assume furthermore that $\left( \mathbf{VK}_{0}\right) 
$, $\left( \mathbf{VK}_{\infty }\right) $ and $\left( \mathbf{f}_{2}\right) $
hold with 
\begin{equation}
\beta _{0},\beta _{\infty }<1,\qquad \alpha _{0}\geq \alpha _{1}^{\left(
0\right) },\qquad \alpha _{\infty }<\alpha _{1}^{\left( \infty \right)
},\qquad \max \left\{ 2\beta _{\infty },q_{0},q_{*}\right\} <q<2,
\label{THM:RNsub1: ineqs}
\end{equation}
where 
\[
\alpha _{1}^{\left( 0\right) }=\alpha _{1}\left( \beta _{0},2\right) ,\quad
\alpha _{1}^{\left( \infty \right) }=\alpha _{1}\left( \beta _{\infty
},2\right) ,\quad q_{0}=q_{0}\left( \alpha _{0},\beta _{0},2\right) \quad 
\text{and}\quad q_{*}=q_{*}\left( \alpha _{\infty },\beta _{\infty
},2\right) .
\]
Then problem $\left( P_{Q}\right) $ has a nonnegative radial solution $u\neq
0$. If $V$ also satisfies $\left( \mathbf{V}_{0}\right) $, then the same
result holds with $\alpha _{1}^{\left( 0\right) }=\alpha _{1}\left( \beta
_{0},\gamma _{0}\right) $ and $q_{0}=q_{0}\left( \alpha _{0},\beta
_{0},\gamma _{0}\right) $, provided that $\alpha _{0}>\alpha _{1}^{\left(
0\right) }$ if $\gamma _{0}\geq N$. If $V$ also satisfies $\left( \mathbf{V}%
_{\infty }\right) $, then we can take $\alpha _{1}^{\left( \infty \right)
}=\alpha _{1}\left( \beta _{\infty },\gamma _{\infty }\right) $ and $%
q_{*}=q_{*}\left( \alpha _{\infty },\beta _{\infty },\gamma _{\infty
}\right) $.
\end{theorem}

\begin{theorem}
\label{THM:RNsub2}Let $\Omega =\mathbb{R}^{N}$ and let $Q\in L^{2}(\mathbb{R}%
_{+},r^{N+1}dr)$, $Q\geq 0$. Assume that $f$ satisfies $\left( \mathbf{F}%
_{3}\right) $, or that $Q$ does not vanish almost everywhere in $\left(
r_{1},r_{2}\right) $. Assume furthermore that $\left( \mathbf{VK}_{0}\right) 
$, $\left( \mathbf{VK}_{\infty }\right) $ and $\left( \mathbf{f}_{1}\right) $
hold with 
\begin{equation}
\beta _{0},\beta _{\infty }<1,\qquad \max \left\{ \alpha _{2},\alpha
_{3}\right\} <\alpha _{0}<\alpha _{1}^{\left( 0\right) },\qquad \alpha
_{\infty }<\alpha _{1}^{\left( \infty \right) },  \label{THM:RNsub2: ineqs0}
\end{equation}
\begin{equation}
\max \left\{ 1,2\beta _{0}\right\} <q_{1}<q_{*}^{\left( 0\right) },\qquad
\max \left\{ 1,2\beta _{\infty },q_{*}^{\left( \infty \right) }\right\}
<q_{2}<2,  \label{THM:RNsub2: ineqs}
\end{equation}
where 
\[
\alpha _{2}=\alpha _{2}\left( \beta _{0}\right) ,\quad \alpha _{3}=\alpha
_{3}\left( \beta _{0},2\right) ,\quad \alpha _{1}^{\left( 0\right) }=\alpha
_{1}\left( \beta _{0},2\right) ,\quad \alpha _{1}^{\left( \infty \right)
}=\alpha _{1}\left( \beta _{\infty },2\right) ,
\]
\[
q_{*}^{\left( 0\right) }=q_{*}\left( \alpha _{0},\beta _{0},2\right) \quad 
\text{and}\quad q_{*}^{\left( \infty \right) }=q_{*}\left( \alpha _{\infty
},\beta _{\infty },2\right) .
\]
If $Q$ does not vanish almost everywhere in $\left( r_{1},r_{2}\right) $, $%
q_{2}=2$ is also allowed in (\ref{THM:RNsub2: ineqs}). Then problem $\left(
P_{Q}\right) $ has a nonnegative radial solution $u\neq 0$. If $V$ also
satisfies $\left( \mathbf{V}_{0}\right) $ with $2<\gamma _{0}<N$, then we
can take $\alpha _{1}^{\left( 0\right) }=\alpha _{1}\left( \beta _{0},\gamma
_{0}\right) $, $\alpha _{3}=\alpha _{3}\left( \beta _{0},\gamma _{0}\right) $
and $q_{*}^{\left( 0\right) }=q_{*}\left( \alpha _{0},\beta _{0},\gamma
_{0}\right) $. If $V$ also satisfies $\left( \mathbf{V}_{\infty }\right) $,
then we can take $\alpha _{1}^{\left( \infty \right) }=\alpha _{1}\left(
\beta _{\infty },\gamma _{\infty }\right) $ and $q_{*}^{\left( \infty
\right) }=q_{*}\left( \alpha _{\infty },\beta _{\infty },\gamma _{\infty
}\right) $.
\end{theorem}

\begin{remark}
\label{RMK:RNsub}\quad 

\begin{enumerate}
\item  The inequalities $\max \left\{ \alpha _{2},\alpha _{3}\right\}
<\alpha _{1}^{\left( 0\right) }$, $\max \left\{ 1,2\beta _{0}\right\}
<q_{*}^{\left( 0\right) }$ and $\max \{1,2\beta _{\infty },q_{*}^{\left(
\infty \right) }\}<2$ in (\ref{THM:RNsub2: ineqs0})-(\ref{THM:RNsub2: ineqs}%
) and $\max \left\{ 2\beta _{\infty },q_{0},q_{*}\right\} <2$ in (\ref
{THM:RNsub1: ineqs}) are ensured by the other hypotheses of Theorems \ref
{THM:RNsub1} and \ref{THM:RNsub2}, so that they are not further assumptions.

\item  \label{RMK:RNsub: useless}As $\left( \mathbf{VK}_{0}\right) $ remains
true if we lower $\alpha _{0}$, the assumption $\alpha _{0}<\alpha
_{1}^{\left( 0\right) }$ ($=\alpha _{1}\left( \beta _{0},\gamma _{0}\right) $%
, $2\leq \gamma _{0}<N$) in (\ref{THM:RNsub2: ineqs0}) is not a restriction.
Nevertheless, if the hypotheses of Theorem \ref{THM:RNsub2} are satisfied
and $\left( \mathbf{VK}_{0}\right) $ holds with $\alpha _{0}\geq \alpha
_{1}^{\left( 0\right) }$, it is never convenient to reduce $\alpha _{0}$ and
apply Theorem \ref{THM:RNsub2}, since one can always apply Theorem \ref
{THM:RNsub1} and get a better result (because $\left( \mathbf{f}_{1}\right) $
with $q_{1},q_{2}$ satisfying (\ref{THM:RNsub2: ineqs}) implies $\left( 
\mathbf{f}_{2}\right) $ for some $q$ satisfying (\ref{THM:RNsub1: ineqs})).
In other words, Theorem \ref{THM:RNsub2} is useful with respect to Theorem 
\ref{THM:RNsub1} only when $\left( \mathbf{VK}_{0}\right) $ does not hold
for some $\alpha _{0}\geq \alpha _{1}^{\left( 0\right) }$.

\item  \label{RMK:RNsub: 1improve}For $\beta \in \left[ 0,1\right) $ fixed, $%
\alpha _{1}\left( \beta ,\gamma \right) $ is continuous and strictly
decreasing in $\gamma \in \mathbb{R}$ and one can check that $q_{0}\left(
\alpha ,\beta ,\gamma \right) $, defined on the set $\left\{ \left( \alpha
,\gamma \right) :\gamma \geq 2,\,\alpha \geq \alpha _{1}\left( \beta ,\gamma
\right) \right\} $, is continuous and decreasing both in $\gamma $ for $%
\alpha $ fixed and in $\alpha $ for $\gamma $ fixed. Similarly, the function
defined on $\left\{ \left( \alpha ,\gamma \right) :\gamma \leq 2,\,\alpha
<\alpha _{1}\left( \beta ,\gamma \right) \right\} $ by $\max \left\{ 2\beta
,q_{*}\left( \alpha ,\beta ,\gamma \right) \right\} $ is increasing and
continuous both in $\gamma $ for $\alpha $ fixed and in $\alpha $ for $%
\gamma $ fixed.

This shows that Theorem \ref{THM:RNsub1} improves under assumption $\left( 
\mathbf{V}_{0}\right) $, or $\left( \mathbf{V}_{\infty }\right) $, or both.

Moreover, as in Remark \ref{RMK:RNsuper}.\ref{RMK:RNsuper: improve}, the
theorem is still true if we replace $\alpha _{0},\alpha _{\infty },\gamma
_{0},\gamma _{\infty }$ with $\overline{\alpha }_{0},\underline{\alpha }%
_{\infty },\overline{\gamma }_{0},\underline{\gamma }_{\infty }$ in $\alpha
_{1}^{\left( 0\right) },\alpha _{1}^{\left( \infty \right) },q_{0},q_{*}$,
and $\alpha _{0}\geq \alpha _{1}^{\left( 0\right) }$ with $\alpha
_{0}>\alpha _{1}^{\left( 0\right) }$ in (\ref{THM:RNsub1: ineqs}).

\item  The same monotonicities in $\gamma $ of Remark \ref{RMK:RNsub}.\ref
{RMK:RNsub: 1improve}, together with the fact that $\max \left\{ \alpha
_{2}\left( \beta \right) ,\alpha _{3}\left( \beta ,\gamma \right) \right\} $
and $q_{*}\left( \alpha ,\beta ,\gamma \right) $ are respectively decreasing
and strictly increasing in $\gamma \in \left[ 2,N\right) $ for $\beta \in
\left[ 0,1\right) $ and $\alpha >\alpha _{2}\left( \beta \right) $ fixed,
show that Theorem \ref{THM:RNsub2} improves under assumption $\left( \mathbf{%
V}_{0}\right) $, or $\left( \mathbf{V}_{\infty }\right) $, or both.
Moreover, as in Remarks \ref{RMK:RNsuper}.\ref{RMK:RNsuper: improve} and \ref
{RMK:RNsub}.\ref{RMK:RNsub: 1improve}, a version of the theorem with $\alpha
_{0},\alpha _{\infty },\gamma _{0},\gamma _{\infty }$ replaced by $\overline{%
\alpha }_{0},\underline{\alpha }_{\infty },\overline{\gamma }_{0},\underline{%
\gamma }_{\infty }$ also holds, the details of which we leave to the
interested reader.

\item  \label{RMK:RNsub: single}In Theorem \ref{THM:RNsub2}, it may happen
that $\max \{1,2\beta _{\infty },q_{*}^{\left( \infty \right)
}\}<q_{*}^{\left( 0\right) }$. In this case, one can take $q_{1}=q_{2}$ and
a single-power growth condition on $f$ is thus enough to apply the theorem
and get existence. By the way, $\alpha _{0}<\alpha _{1}^{\left( 0\right) }$
and $\gamma _{0}<N$ imply $q_{*}^{\left( 0\right) }<2$ (and therefore $%
q_{1}<2$), so that the linear case $q_{1}=q_{2}=2$ is always excluded.

\item  In Theorems \ref{THM:RNsub1} and \ref{THM:RNsub2}, the requirement $%
Q\in L^{2}(\mathbb{R}_{+},r^{N+1}dr)$ just plays the role of ensuring that the
linear operator $u\mapsto \int_{\mathbb{R}^{N}}Q\left( \left| x\right| \right)
u\,dx$ is continuous on $H_{0,V}^{1}(\mathbb{R}^{N})$ (see Section \ref{SEC: pfs} below).
Therefore it can be replaced by any other condition giving the same
property, e.g., $Q\in L^{2N/\left( N+2\right) }(\mathbb{R}_{+},r^{N-1}dr)$
or $QV^{-1/2}\in L^{2}(\mathbb{R}_{+},r^{N-1}dr)$.

\item  According to the proofs, the solution $u$ of Theorems \ref{THM:RNsub1}
and \ref{THM:RNsub2} satisfies $%
\displaystyle%
I\left( u\right) =\min_{v\in H_{0,V,\mathrm{r}}^{1}(\mathbb{R}^{N}),\,v\geq
0}I\left( v\right) $, where 
\begin{equation}
I\left( v\right) :=\frac{1}{2}\int_{\mathbb{R}^{N}}\left( \left| \nabla
v\right| ^{2}+V\left( \left| x\right| \right) v^{2}\right) dx-\int_{\mathbb{R}%
^{N}}\left( K\left( \left| x\right| \right) F\left( v\right) +Q\left( \left|
x\right| \right) v\right) dx.  \label{I(v):=}
\end{equation}
Moreover, if $Q$ does not vanish almost everywhere in $\left(
r_{1},r_{2}\right) $, both the theorems still work even without assuming $%
Q\geq 0$ (use Theorem \ref{THM:ex-sub} instead of Corollary \ref{COR:ex-sub}
in the proof), but we cannot ensure anymore that $u$ is nonnegative. In this
case, the solution satisfies $%
\displaystyle%
I\left( u\right) =\min_{v\in H_{0,V,\mathrm{r}}^{1}(\mathbb{R}^{N})}I\left(
v\right) $.
\end{enumerate}
\end{remark}

Exactly as in the superlinear case, the same arguments leading to Theorems 
\ref{THM:RNsub1} and \ref{THM:RNsub2} also yield existence results for the
Dirichlet problem in bounded balls or exterior radial domains, where,
respectively, only assumptions on $V$ and $K$ near the origin or at infinity
are needed. In both cases, a single-power growth condition on the
nonlinearity is sufficient. The precise statements are left to the
interested reader.\medskip 

We conclude with a multiplicity result for problem $\left( P_{0}\right) $,
which, in the superlinear case, requires the following assumption,
complementary to $\left( \mathbf{f}_{1}\right) $:

\begin{itemize}
\item[$\left( \mathbf{f}_{1}^{\prime }\right) $]  $\exists q_{1},q_{2}>1$
such that 
\[
\inf_{t>0}\,\frac{f\left( t\right) }{\min \left\{
t^{q_{1}-1},t^{q_{2}-1}\right\} }>0.
\]
\end{itemize}

\begin{theorem}
\label{THM:multiple}\emph{(i)} Under the same assumptions of each of
Theorems \ref{THM:RNsuper}, \ref{THM:Bsuper} and \ref{THM:BCsuper}, if $f$
is also odd and satisfies $\left( \mathbf{f}_{1}^{\prime }\right) $ (with
the same exponents $q_{1},q_{2}$ of $\left( \mathbf{f}_{1}\right) $), then
problem $\left( P_{0}\right) $ has infinitely many radial solutions. \emph{%
(ii)} Under the same assumptions of each of Theorems \ref{THM:RNsub1} and 
\ref{THM:RNsub2} with $Q=0$, if $f$ is also odd, then problem $\left(
P_{0}\right) $ has infinitely many radial solutions.
\end{theorem}

\begin{remark}
The infinitely many solutions of Theorem \ref{THM:multiple} form a sequence $%
\left\{ u_{n}\right\} $ such that $I\left( u_{n}\right) \rightarrow +\infty $
in the superlinear case and $I\left( u_{n}\right) \rightarrow 0$ in the
sublinear one, where $I$ is the functional defined in (\ref{I(v):=}).
\end{remark}

In \cite[Section 3]{BGR-p1} and \cite{GR15}, we have discussed many examples
of pairs of potentials $V,K$ and nonlinearities $f$ satifying our
hypotheses. In the same papers, we have also compared such hypotheses with
the assumptions of some of the main related results in the previous
literature, showing essentially that conditions $\left( \mathbf{VK}%
_{0}\right) $, $\left( \mathbf{VK}_{\infty }\right) $ and $\left( \mathbf{f}%
_{1}\right) $ allow to deal with potentials exhibiting behaviours at zero
and at infinity which are new in the literature (see \cite{BGR-p1} and
Section \ref{SEC: ex} below) and do not need to be compatible with each
other (see both \cite{BGR-p1,GR15}). The same examples, and the same
discussion, can be repetead here, covering many cases not included in
previous papers. In particular, our results for problem $\left( P_{0}\right) 
$ in $\mathbb{R}^{N}$ contain and extend in different directions the results of 
\cite{SuTian12,Su-Wang-Will-p} (for $p=2$) and are complementary to the ones
of \cite{BonMerc11,Su12,Yang-Li-14}. To the best of our knowledge,
assumptions $\left( \mathbf{VK}_{0}\right) $ and $\left( \mathbf{VK}_{\infty
}\right) $ are also new in the study of problem $\left( P_{Q}\right) $ in
bounded or exterior domains.
\smallskip 

The paper is organized as follows. In Section \ref{SEC: ex} we give an
example, complementary to the ones of \cite[Section 3]{BGR-p1}, of
potentials satisfying our hypotheses but not included in the results known
in the literature up to now. In Section \ref{SEC: gen-results}, we introduce
our variational approach to problem $\left( P\right) $ and give some
existence and multiplicity results (Theorems \ref{THM:ex}-\ref{THM:mult-sub}
and Corollary \ref{COR:ex-sub}), which are more general than the ones stated
in this introduction but rely on a less explicit assumption on the
potentials (condition $\left( \mathcal{S}_{q_{1},q_{2}}^{\prime \prime
}\right) $). We also give a suitable symmetric criticality type principle
(Proposition \ref{PROP:symm-crit}), since the Palais' classical one \cite
{Palais} does not apply in this case. Sections \ref{SEC:pf-general} and \ref
{SEC: pfs} are respectively devoted to the proof of the results stated in
Section \ref{SEC: gen-results} and in the Introduction. In particular, the
latter follow from the former, by applying the compactness theorems of \cite
{BGR-p1} (see Lemma \ref{LEM:compactness} below). In the Appendix, we prove
some pointwise estimates for radial Sobolev functions, already used in \cite
{BGR-p1}.


\section{An example \label{SEC: ex}}

A particular feature of our results
is that we do not
necessarily assume hypotheses on $V$ and $K$ separately, but rather on their
ratio: the ratio must have a power-like behaviour at zero and at infinity,
but $V$ and $K$ are not obliged to have such a behaviour. This allows us to
deal with potentials $V,K$ that grow (or vanish) very fast at zero or
infinity, in such a way that they escape the results of the previous
literature but the ratio $K/V$ satisfies our hypotheses.

In particular we can also treat some examples of potentials not included
among those considered in \cite{BonMerc11}, a paper that deals with a very
general class of potentials, the so called \emph{Hardy-Dieudonn\'{e} class},
which also includes the potentials treated in \cite{Su12,SuTian12,Su-Wang-Will-p,Yang-Li-14}. 
In this class, the functions with the fastest growth at infinity are the $n$ times compositions of the
exponential map with itself. Accordingly, let us denote by 
$e_{n}:\left[0,+\infty \right) \rightarrow \mathbb{R}$ 
the function obtained by composing the exponential map with itself $n\geq 0$ times. 
The Hardy-Dieudonn\'{e} class contains $e_{n}$ for all $n$ and these are the mappings with the
fastest growth in the class, so that any function growing faster than every 
$e_{n}$ is not in that class. Let us then define 
$\alpha :\left[ 0,+\infty \right) \rightarrow \mathbb{R}$ by setting 
$\alpha \left( r\right):=e_{n}\left( n\right) $ if $n\leq r<n+1$. 
It is clear that 
\[
\lim_{r\rightarrow +\infty }\frac{\alpha \left( r\right) }{e_{n}\left(
r\right) }=+\infty \quad \textrm{for all }n
\]
and therefore $\alpha $ does not belong to the Hardy-Dieudonn\'{e} class.
Hence, defining 
\[
K\left( r\right) :=\alpha \left( r\right) \quad \text{and}\quad V\left(
r\right) :=\alpha \left( r\right) V_{1}\left( r\right) 
\]
where $V_{1}$ is any potential satisfying $\left( \mathbf{V}\right) $ and
having suitable power-like behavior at zero and infinity, we get a pair of
potentials $V,K$ which satisfy our hypotheses but not those of \cite
{BonMerc11}. Of course, one can build similar examples of potential pairs
which vanish so fast at infinity (or grow or vanish so fast at zero) that
they fall out of the Hardy-Dieudonn\'{e} class, yet their ratio exhibits a
power-like behaviour.

\section{Variational approach and general results \label{SEC: gen-results}}

Let $N\geq 3$ and let $V:\mathbb{R}_{+}\rightarrow \left[ 0,+\infty \right] $
be a measurable function satisfying the following hypothesis:

\begin{itemize}
\item[$\left( h_{0}\right) $]  $V\in L^{1}\left( \left( r_{1},r_{2}\right)
\right) $ for some $r_{2}>r_{1}>0$.
\end{itemize}

\noindent Define the Hilbert spaces 
\begin{eqnarray}
&&H_{V}^{1} 
\begin{array}{c}
:=
\end{array}
H_{V}^{1}\left( \mathbb{R}^{N}\right) 
\begin{array}{c}
:=
\end{array}
\left\{ u\in D^{1,2}\left( \mathbb{R}^{N}\right) :\int_{\mathbb{R}^{N}}V\left(
\left| x\right| \right) u^{2}dx<\infty \right\} ,  \label{H^1_V :=} \\
&&H_{V,\mathrm{r}}^{1} 
\begin{array}{c}
:=
\end{array}
H_{V,\mathrm{r}}^{1}\left( \mathbb{R}^{N}\right) 
\begin{array}{c}
:=
\end{array}
\left\{ u\in H_{V}^{1}\left( \mathbb{R}^{N}\right) :u\left( x\right) =u\left(
\left| x\right| \right) \right\}  \label{H^1_Vr :=}
\end{eqnarray}
where $ D^{1,2}\left( \mathbb{R}^{N}\right)$ is the usual Sobolev space. $H_{V}^{1}$ and $H_{V,\mathrm{r}}^{1}$ are endowed with the following inner product and related norm: 
\begin{equation}
\left( u\mid v\right) :=\int_{\mathbb{R}^{N}}\nabla u\cdot \nabla v\,dx+\int_{%
\mathbb{R}^{N}}V\left( \left| x\right| \right) uv\,dx,\qquad \left\| u\right\|
:=\left( \int_{\mathbb{R}^{N}}\left| \nabla u\right| ^{2}dx+\int_{\mathbb{R}%
^{N}}V\left( \left| x\right| \right) u^{2}dx\right) ^{1/2}.  \label{h struct}
\end{equation}
Of course, $u\left( x\right) =u\left( \left| x\right| \right) $ means that $%
u $ is invariant under the action on $H_{V}^{1}$ of the orthogonal group of $%
\mathbb{R}^{N}$. Note that $H_{V}^{1}$ and $H_{V,\mathrm{r}}^{1}$ are nonzero
by $\left( h_{0}\right) $ and $H_{V}^{1}$ is the space $H_{0,V}^{1}(\mathbb{R}%
^{N})$ defined in the Introduction.

Let $g:\mathbb{R}_{+}\times \mathbb{R}\rightarrow \mathbb{R}$ be a Carath\'{e}odory
function and assume once and for all that there exist $f\in C\left( \mathbb{R};%
\mathbb{R}\right) $ and a measurable function $K:\mathbb{R}_{+}\rightarrow \mathbb{R}%
_{+}$ such that:

\begin{itemize}
\item[$\left( h_{1}\right) $]  $\left| g\left( r,t\right) -g\left(
r,0\right) \right| \leq K\left( r\right) \left| f\left( t\right) \right| $
for almost every $r>0$ and all $t\in \mathbb{R};$

\item[$\left( h_{2}\right) $]  $K\in L_{\mathrm{loc}}^{s}\left( \left(
0,+\infty \right) \right) $ for some $s>\frac{2N}{N+2}.$
\end{itemize}

\noindent Assume furthermore that:

\begin{itemize}
\item[$\left( h_{3}\right) $]  the linear operator $u\mapsto \int_{\mathbb{R}%
^{N}}g\left( \left| x\right| ,0\right) u\,dx$ is continuous on $H_{V}^{1}$
\end{itemize}

\noindent (see also Remark \ref{RMK: h3}). Of course $\left( h_{3}\right) $
will be relevant only if $g\left( \cdot ,0\right) \neq 0$ (meaning that $%
g\left( \cdot ,0\right) $ does not vanish almost everywhere).\smallskip

Define the following functions of $R>0$ and $q>1$: 
\[
\mathcal{S}_{0}\left( q,R\right) := 
\sup_{u\in H_{V,\mathrm{r}}^{1},\,
\left\| u\right\| =1}\,
\int_{B_{R}}K\left( \left| x\right| \right)\left| u\right| ^{q}dx,\quad 
\mathcal{S}_{\infty }\left( q,R\right) := 
\sup_{u\in H_{V,\mathrm{r}}^{1},\,
\left\| u\right\| =1}\,
\int_{\mathbb{R}^{N}\setminus B_{R}}K\left( \left| x\right| \right) \left| u\right| ^{q}dx, 
\]
\begin{eqnarray*}
\mathcal{R}_{0}\left( q,R\right)&:=&
\sup_{u\in H_{V,\mathrm{r}}^{1},\,h\in H_{V}^{1},\,
\left\| u\right\| =\left\| h\right\| =1}\,
\int_{B_{R}}K\left( \left| x\right| \right) \left| u\right| ^{q-1}\left|
h\right| dx, \\
\mathcal{R}_{\infty }\left( q,R\right) &:=&
\sup_{u\in H_{V,\mathrm{r}}^{1},\,h\in H_{V}^{1},\,
\left\| u\right\| =\left\| h\right\| =1}\,
\int_{\mathbb{R}^{N}\setminus B_{R}}K\left( \left| x\right| \right) \left|
u\right| ^{q-1}\left| h\right| dx.
\end{eqnarray*}
Note that $\mathcal{S}_{0}\left( q,\cdot \right) $ and $\mathcal{R}%
_{0}\left( q,\cdot \right) $ are increasing, $\mathcal{S}_{\infty }\left(
q,\cdot \right) $ and $\mathcal{R}_{\infty }\left( q,\cdot \right) $ are
decreasing and all can be infinite at some $R$. Moreover, for every $\left(
q,R\right) $ one has $\mathcal{S}_{0}\left( q,R\right) \leq \mathcal{R}%
_{0}\left( q,R\right) $ and $\mathcal{S}_{\infty }\left( q,R\right) \leq 
\mathcal{R}_{\infty }\left( q,R\right) $.

On the functions $\mathcal{S},\mathcal{R}$ and $f$, we will require suitable
combinations of the following conditions (see also Remarks \ref{RMK:thm:ex}.%
\ref{RMK:thm:ex-2} and \ref{RMK:ex-sub}), where $q_{1},q_{2}$ will be
specified each time:

\begin{itemize}
\item[$\left( f_{q_{1},q_{2}}\right) $]  $\exists M>0$ such that $\left|
f\left( t\right) \right| \leq M\min \left\{ \left| t\right|
^{q_{1}-1},\left| t\right| ^{q_{2}-1}\right\} $ for all $t\in \mathbb{R};$

\item[$\left( \mathcal{S}_{q_{1},q_{2}}^{\prime }\right) $]  $\exists
R_{1},R_{2}>0$ such that $\mathcal{S}_{0}\left( q_{1},R_{1}\right) <\infty $
and $\mathcal{S}_{\infty }\left( q_{2},R_{2}\right) <\infty ;$

\item[$\left( \mathcal{S}_{q_{1},q_{2}}^{\prime \prime }\right) $]  $%
\displaystyle%
\lim_{R\rightarrow 0^{+}}\mathcal{S}_{0}\left( q_{1},R\right)
=\lim_{R\rightarrow +\infty }\mathcal{S}_{\infty }\left( q_{2},R\right) =0;$

\item[$\left( \mathcal{R}_{q_{1},q_{2}}\right) $]  $\exists R_{1},R_{2}>0$
such that $\mathcal{R}_{0}\left( q_{1},R_{1}\right) <\infty $ and $\mathcal{R%
}_{\infty }\left( q_{2},R_{2}\right) <\infty .$
\end{itemize}

We set $G\left( r,t\right) :=\int_{0}^{t}g\left( r,s\right) ds$ and 
\begin{equation}
I\left( u\right) :=\frac{1}{2}\left\| u\right\| ^{2}-\int_{\mathbb{R}%
^{N}}G\left( \left| x\right| ,u\right) dx.  \label{I:=}
\end{equation}
From the embedding results of \cite{BGR-p1} and the results of \cite{BPR}
about Nemytski\u{\i} operators on the sum of Lebesgue spaces (see Section 
\ref{SEC:pf-general} below for some recallings on such spaces), we get the
following differentiability result.

\begin{proposition}
\label{PROP:diff}Assume that there exist $q_{1},q_{2}>1$ such that $\left(
f_{q_{1},q_{2}}\right) $ and $\left( \mathcal{S}_{q_{1},q_{2}}^{\prime
}\right) $ hold. Then (\ref{I:=}) defines a $C^{1}$ functional on $H_{V,%
\mathrm{r}}^{1}$, with Fr\'{e}chet derivative at any $u\in H_{V,\mathrm{r}%
}^{1}$ given by 
\begin{equation}
I^{\prime }\left( u\right) h=\int_{\mathbb{R}^{N}}\left( \nabla u\cdot \nabla
h+V\left( \left| x\right| \right) uh\right) dx-\int_{\mathbb{R}^{N}}g\left(
\left| x\right| ,u\right) h\,dx,\quad \forall h\in H_{V,\mathrm{r}}^{1}.
\label{PROP:diff: I'(u)h=}
\end{equation}
\end{proposition}

Proposition \ref{PROP:diff} ensures that the critical points of $I:H_{V,%
\mathrm{r}}^{1}\rightarrow \mathbb{R}$ satisfy (\ref{weak solution}) (with $%
\Omega =\mathbb{R}^{N}$) for all $h\in H_{V,\mathrm{r}}^{1}$. The next result
shows that such critical points are actually weak solutions to problem $%
\left( P\right) $ (with $\Omega =\mathbb{R}^{N}$), provided that the slightly
stronger version $\left( \mathcal{R}_{q_{1},q_{2}}\right) $ of condition $%
\left( \mathcal{S}_{q_{1},q_{2}}^{\prime }\right) $ holds. Observe that the
classical Palais' Principle of Symmetric Criticality \cite{Palais} does not
apply in this case, because we do not know whether or not $I$ is
differentiable, not even well defined, on the whole space $H_{V}^{1}$.

\begin{proposition}
\label{PROP:symm-crit}Assume that there exist $q_{1},q_{2}>1$ such that $%
\left( f_{q_{1},q_{2}}\right) $ and $\left( \mathcal{R}_{q_{1},q_{2}}\right) 
$ hold. Then every critical point $u$ of $I:H_{V,\mathrm{r}}^{1}\rightarrow 
\mathbb{R}$ satisfies 
\begin{equation}
\int_{\mathbb{R}^{N}}\nabla u\cdot \nabla h\,dx+\int_{\mathbb{R}^{N}}V\left(
\left| x\right| \right) uh\,dx=\int_{\mathbb{R}^{N}}g\left( \left| x\right|
,u\right) h\,dx,\quad \forall h\in H_{V}^{1}  \label{PROP:symm-crit: th}
\end{equation}
(i.e., $u$ is a weak solution to problem $\left( P\right) $ with $\Omega =%
\mathbb{R}^{N}$).
\end{proposition}

By Proposition \ref{PROP:symm-crit}, the problem of radial weak solutions to 
$\left( P\right) $ (with $\Omega =\mathbb{R}^{N}$) reduces to the study of
critical points of $I:H_{V,\mathrm{r}}^{1}\rightarrow \mathbb{R}$. Concerning
the case of superlinear nonlinearities, we have the following existence and
multipilicity results.

\begin{theorem}
\label{THM:ex}Assume $g\left( \cdot ,0\right) =0$ and assume that there
exist $q_{1},q_{2}>2$ such that $\left( f_{q_{1},q_{2}}\right) $ and $\left( 
\mathcal{S}_{q_{1},q_{2}}^{\prime \prime }\right) $ hold. Assume furthermore
that $g$ satisfies:

\begin{itemize}
\item[$\left( g_{1}\right) $]  $\exists \theta >2$ such that $0\leq \theta
G\left( r,t\right) \leq g\left( r,t\right) t$ for almost every $r>0$ and all 
$t\geq 0;$

\item[$\left( g_{2}\right) $]  $\exists t_{0}>0$ such that $G\left(
r,t_{0}\right) >0$ for almost every $r>0.$
\end{itemize}

\noindent If $K\left( \left| \cdot \right| \right) \in L^{1}(\mathbb{R}^{N})$,
we can replace assumptions $\left( g_{1}\right) $-$\left( g_{2}\right) $
with:

\begin{itemize}
\item[$\left( g_{3}\right) $]  $\exists \theta >2$ and $\exists t_{0}>0$
such that $0<\theta G\left( r,t\right) \leq g\left( r,t\right) t$ for almost
every $r>0$ and all $t\geq t_{0}.$
\end{itemize}

\noindent Then the functional $I:H_{V,\mathrm{r}}^{1}\rightarrow \mathbb{R}$
has a nonnegative critical point $u\neq 0$.
\end{theorem}

\begin{remark}
\label{RMK:thm:ex}\quad 

\begin{enumerate}
\item  \label{RMK:thm:ex-1}Assumptions $\left( g_{1}\right) $ and $\left(
g_{2}\right) $ imply $\left( g_{3}\right) $, so that, in Theorem \ref{THM:ex}%
, the information $K\left( \left| \cdot \right| \right) \in L^{1}(\mathbb{R}%
^{N})$ actually allows weaker hypotheses on the nonlinearity.

\item  \label{RMK:thm:ex-2}In Theorem \ref{THM:ex}, assumptions $\left(
h_{1}\right) $ and $\left( f_{q_{1},q_{2}}\right) $ need only to hold for $%
t\geq 0$. Indeed, all the hypotheses of the theorem still hold true if we
replace $g\left( r,t\right) $ with $\chi _{\mathbb{R}_{+}}\left( t\right)
g\left( r,t\right) $ ($\chi _{\mathbb{R}_{+}}$ is the characteristic function
of $\mathbb{R}_{+}$) and this can be done without restriction since the theorem
concerns \emph{nonnegative} critical points.
\end{enumerate}
\end{remark}

\begin{theorem}
\label{THM:mult}Assume that there exist $q_{1},q_{2}>2$ such that $\left(
f_{q_{1},q_{2}}\right) $ and $\left( \mathcal{S}_{q_{1},q_{2}}^{\prime
\prime }\right) $ hold. Assume furthermore that:

\begin{itemize}
\item[$\left( g_{4}\right) $]  $\exists m>0$ such that $G\left( r,t\right)
\geq mK\left( r\right) \min \left\{ t^{q_{1}},t^{q_{2}}\right\} $ for almost
every $r>0$ and all $t\geq 0;$

\item[$\left( g_{5}\right) $]  $g\left( r,t\right) =-g\left( r,-t\right) $
for almost every $r>0$ and all $t\geq 0.$
\end{itemize}

\noindent Finally, assume that $g$ satisfies $\left( g_{1}\right) $, or that 
$K\left( \left| \cdot \right| \right) \in L^{1}(\mathbb{R}^{N})$ and $g$
satisfies $\left( g_{3}\right) $. Then the functional $I:H_{V,\mathrm{r}%
}^{1}\rightarrow \mathbb{R}$ has a sequence of critical points $\left\{
u_{n}\right\} $ such that $I\left( u_{n}\right) \rightarrow +\infty $.
\end{theorem}

\begin{remark}
The condition $g\left( \cdot ,0\right) =0$ is implicit in Theorem \ref
{THM:mult} (and in Theorem \ref{THM:mult-sub} below), as it follows from
assumption $\left( g_{5}\right) $.
\end{remark}

As to sublinear nonlinearities, we will prove the following results.

\begin{theorem}
\label{THM:ex-sub}Assume that there exist $q_{1},q_{2}\in \left( 1,2\right) $
such that $\left( f_{q_{1},q_{2}}\right) $ and $\left( \mathcal{S}%
_{q_{1},q_{2}}^{\prime \prime }\right) $ hold. Assume furthermore that $g$
satisfies at least one of the following conditions:

\begin{itemize}
\item[$\left( g_{6}\right) $]  $\exists \theta <2$ and $\exists t_{0},m>0$
such that $G\left( r,t\right) \geq mK\left( r\right) t^{\theta }$ for almost
every $r>0$ and all $0\leq t\leq t_{0};$

\item[$\left( g_{7}\right) $]  $g\left( \cdot ,0\right) $ does not vanish
almost everywhere in $\left( r_{1},r_{2}\right) .$
\end{itemize}

\noindent If $\left( g_{7}\right) $ holds, we also allow the case $\max
\left\{ q_{1},q_{2}\right\} =2>\min \left\{ q_{1},q_{2}\right\} >1$. Then
there exists $u\neq 0$ such that 
\[
I\left( u\right) =\min_{v\in H_{V,\mathrm{r}}^{1}}I\left( v\right) .
\]
\end{theorem}

If $g\left( \cdot ,t\right) \geq 0$ almost everywhere for all $t<0$, the
minimizer $u$ of Theorem \ref{THM:ex-sub} is nonnegative, since a standard
argument shows that all the critical points of $I$ are nonnegative (test $%
I^{\prime }\left( u\right) $ with the negative part $u_{-}$ and get $%
I^{\prime }\left( u\right) u_{-}=-\left\| u_{-}\right\| ^{2}=0$). The next
corollary gives a nonnegative critical point just asking $g\left(\cdot ,0\right) \geq 0$.

\begin{corollary}
\label{COR:ex-sub}Assume the same hypotheses of Theorem \ref{THM:ex-sub}. If 
$g\left( \cdot ,0\right) \geq 0$ almost everywhere, then $I:H_{V,\mathrm{r}%
}^{1}\rightarrow \mathbb{R}$ has a nonnegative critical point $\widetilde{u}%
\neq 0$ satisfying 
\begin{equation}
I\left( \widetilde{u}\right) =\min_{u\in H_{V,\mathrm{r}}^{1},\,u\geq
0}I\left( u\right) .  \label{COR:ex-sub: th}
\end{equation}
\end{corollary}

\begin{remark}
\label{RMK:ex-sub}\quad 

\begin{enumerate}
\item  In Theorem \ref{THM:ex-sub} and Corollary \ref{COR:ex-sub}, the case $%
\max \left\{ q_{1},q_{2}\right\} =2>\min \left\{ q_{1},q_{2}\right\} >1$
cannot be considered under assumption $\left( g_{6}\right) $, since $\left(
g_{6}\right) $ and $\left( f_{q_{1},q_{2}}\right) $ imply $\max \left\{
q_{1},q_{2}\right\} \leq \theta <2$.

\item  \label{RMK:ex-sub-2}Checking the proof, one sees that Corollary \ref
{COR:ex-sub} actually requires that assumptions $\left( h_{1}\right) $ and $%
\left( f_{q_{1},q_{2}}\right) $ hold only for $t\geq 0$, which is consistent
with the concern of the result about \emph{nonnegative} critical points.
\end{enumerate}
\end{remark}

\begin{theorem}
\label{THM:mult-sub}Assume that there exist $q_{1},q_{2}\in \left(
1,2\right) $ such that $\left( f_{q_{1},q_{2}}\right) $ and $\left( \mathcal{%
S}_{q_{1},q_{2}}^{\prime \prime }\right) $ hold. Assume furthermore that $g$
satisfies $\left( g_{5}\right) $ and $\left( g_{6}\right) $. Then the
functional $I:H_{V,\mathrm{r}}^{1}\rightarrow \mathbb{R}$ has a sequence of
critical points $\left\{ u_{n}\right\} $ such that $I\left( u_{n}\right) <0$
and $I\left( u_{n}\right) \rightarrow 0$.
\end{theorem}

\begin{remark}
\label{RMK: h3}If the linear operator of assumption $\left( h_{3}\right) $
is just continuous on $H_{V,\mathrm{r}}^{1}$, then Proposition \ref
{PROP:symm-crit} fails, but all the other results of this section remain
valid (as can be easily seen by checking the proofs). This is especially
relevant in connection with the radial estimates satisfied by the $H_{V,%
\mathrm{r}}^{1}$ mappings (see Appendix), which ensure that $\left(
h_{3}\right) $ holds on $H_{V,\mathrm{r}}^{1}$ provided that $g\left( \left|
\cdot \right| ,0\right) $ belongs to $L_{\mathrm{loc}}^{1}(\mathbb{R}%
^{N}\setminus \left\{ 0\right\} )$ and satisfies suitable decay (or growth)
conditions at zero and at infinity.
\end{remark}

\section{Proof of the general results \label{SEC:pf-general}}

In this section we keep the notation and assumptions of the preceding
section. Denoting $L_{K}^{p}\left( E\right) :=L^{p}\left( E,K\left( \left|
x\right| \right) dx\right) $ for any measurable set $E\subseteq \mathbb{R}^{N}$%
, we will make frequent use of the sum space 
\[
L_{K}^{p_{1}}+L_{K}^{p_{2}}:=\left\{ u_{1}+u_{2}:u_{1}\in
L_{K}^{p_{1}}\left( \mathbb{R}^{N}\right) ,\,u_{2}\in L_{K}^{p_{2}}\left( \mathbb{R%
}^{N}\right) \right\} ,\quad 1<p_{i}<\infty . 
\]
We recall from \cite{BPR} that such a space can be characterized as the set
of measurable mappings $u:\mathbb{R}^{N}\rightarrow \mathbb{R}$ for which there
exists a measurable set $E\subseteq \mathbb{R}^{N}$ such that $u\in
L_{K}^{p_{1}}\left( E\right) \cap L_{K}^{p_{2}}\left( E^{c}\right) $. It is
a Banach space with respect to the norm 
\[
\left\| u\right\| _{L_{K}^{p_{1}}+L_{K}^{p_{2}}}:=\inf_{u_{1}+u_{2}=u}\max
\left\{ \left\| u_{1}\right\| _{L_{K}^{p_{1}}(\mathbb{R}^{N})},\left\|
u_{2}\right\| _{L_{K}^{p_{2}}(\mathbb{R}^{N})}\right\} 
\]
and the continuous embedding $L_{K}^{p}\hookrightarrow
L_{K}^{p_{1}}+L_{K}^{p_{2}}$ holds for all $p\in \left[ \min \left\{
p_{1},p_{2}\right\} ,\max \left\{ p_{1},p_{2}\right\} \right] $. Moreover,
for every $u\in L_{K}^{p_{1}}+L_{K}^{p_{2}}$ one has 
\begin{equation}
\left\| u\right\| _{L_{K}^{p_{1}}+L_{K}^{p_{2}}}\leq \left\| u\right\|
_{L_{K}^{\min \left\{ p_{1},p_{2}\right\} }(\Lambda _{u})}+\left\| u\right\|
_{L_{K}^{\max \left\{ p_{1},p_{2}\right\} }(\Lambda _{u}^{c})},\quad \text{%
where}\quad \Lambda _{u}:=\left\{ x\in \mathbb{R}^{N}:\left| u\left( x\right)
\right| >1\right\}  \label{Vu}
\end{equation}
(see \cite[Corollary 2.19]{BPR}).\bigskip

\noindent \textbf{Proof of Proposition \ref{PROP:diff}.}\quad 
On the one hand, by $\left( \mathcal{S}_{q_{1},q_{2}}^{\prime }\right) $ and 
\cite[Theorem 1]{BGR-p1}, the embedding $H_{V,\mathrm{r}}^{1}\hookrightarrow
L_{K}^{q_{1}}+L_{K}^{q_{2}}$ is continuous. On the other hand, by $\left(
h_{1}\right) $, $\left( f_{q_{1},q_{2}}\right) $ and \cite[Proposition 3.8]
{BPR}, the functional 
\begin{equation}
u\mapsto \int_{\mathbb{R}^{N}}\left( G\left( \left| x\right| ,u\right) -g\left(
\left| x\right| ,0\right) u\right) dx  \label{NemOp}
\end{equation}
is of class $C^{1}$ on $L_{K}^{q_{1}}+L_{K}^{q_{2}}$ and its Fr\'{e}chet
derivative at any $u$ is given by 
\[
h\in L_{K}^{q_{1}}+L_{K}^{q_{2}}\mapsto \int_{\mathbb{R}^{N}}\left( g\left(
\left| x\right| ,u\right) -g\left( \left| x\right| ,0\right) \right) h\,dx. 
\]
Hence, by $\left( h_{3}\right) $, we conclude that $I\in C^{1}(H_{V,\mathrm{r%
}}^{1})$ and that (\ref{PROP:diff: I'(u)h=}) holds.%
\endproof
\bigskip

\noindent \textbf{Proof of Proposition \ref{PROP:symm-crit}.}\quad 
Let $u\in H_{V,\mathrm{r}}^{1}$. By the monotonicity of $\mathcal{R}_{0}$ and $%
\mathcal{R}_{\infty }$, it is not restrictive to assume $R_{1}<R_{2}$ in
hypothesis $\left( \mathcal{R}_{q_{1},q_{2}}\right) $. So, by \cite[Lemma 1]
{BGR-p1}, there exists a constant $C>0$ (dependent on $u$) such that for all 
$h\in H_{V}^{1}$ we have 
\[
\int_{B_{R_{2}}\setminus B_{R_{1}}}K\left( \left| x\right| \right) \left|
u\right| ^{q_{1}-1}\left| h\right| dx\leq C\left\| h\right\| 
\]
and therefore, by $\left( h_{1}\right) $ and $\left( f_{q_{1},q_{2}}\right) $%
, 
\begin{eqnarray*}
\int_{\mathbb{R}^{N}}\left| \,g\left( \left| x\right| ,u\right) -g\left( \left|
x\right| ,0\right) \,\right| \left| h\right| dx &\leq &\int_{\mathbb{R}%
^{N}}K\left( \left| x\right| \right) \left| f\left( u\right) \right| \left|
h\right| dx\leq M\int_{\mathbb{R}^{N}}K\left( \left| x\right| \right) \min
\{\left| u\right| ^{q_{1}-1},\left| u\right| ^{q_{2}-1}\}\left| h\right| dx
\\
&\leq &M\left( \int_{B_{R_{1}}}K\left( \left| x\right| \right) \left|
u\right| ^{q_{1}-1}\left| h\right| dx+\int_{B_{R_{2}}^{c}}K\left( \left|
x\right| \right) \left| u\right| ^{q_{2}-1}\left| h\right| dx\right. + \\
&&+\left. \int_{B_{R_{2}}\setminus B_{R_{1}}}K\left( \left| x\right| \right)
\left| u\right| ^{q_{1}-1}\left| h\right| dx\right) \\
&\leq &M\left( \left\| u\right\| ^{q_{1}-1}\left\| h\right\|
\int_{B_{R_{1}}}K\left( \left| x\right| \right) \frac{\left| u\right|
^{q_{1}-1}}{\left\| u\right\| ^{q_{1}-1}}\frac{\left| h\right| }{\left\|
h\right\| }dx\right. + \\
&&+\left. \left\| u\right\| ^{q_{2}-1}\left\| h\right\|
\int_{B_{R_{2}}^{c}}K\left( \left| x\right| \right) \frac{\left| u\right|
^{q_{2}-1}}{\left\| u\right\| ^{q_{2}-1}}\frac{\left| h\right| }{\left\|
h\right\| }dx+C\left\| h\right\| \right) \\
&\leq &M\left( \left\| u\right\| ^{q_{1}-1}\mathcal{R}_{0}\left(
q_{1},R_{1}\right) +\left\| u\right\| ^{q_{2}-1}\mathcal{R}_{\infty }\left(
q_{2},R_{2}\right) +C\right) \left\| h\right\| .
\end{eqnarray*}
Together with $\left( h_{3}\right) $, this gives that the linear operator 
\[
T\left( u\right) h:=\int_{\mathbb{R}^{N}}\left( \nabla u\cdot \nabla h+V\left(
\left| x\right| \right) uh\right) dx-\int_{\mathbb{R}^{N}}g\left( \left|
x\right| ,u\right) h\,dx 
\]
is well defined and continuous on $H_{V}^{1}$. Hence, by Riesz
representation theorem, there exists a unique $\tilde{u}\in H_{V}^{1}$ such
that $T\left( u\right) h=\left( \tilde{u}\mid h\right) $ for all $h\in
H_{V}^{1}$, where $\left( \cdot \mid \cdot \right) $ is the inner product of 
$H_{V}^{1}$ defined in (\ref{h struct}). Denoting by $O\left( N\right) $ the
orthogonal group of $\mathbb{R}^{N}$, by means of obvious changes of variables
it is easy to see that for every $h\in H_{V}^{1}$ and $g\in O\left( N\right) 
$ one has $\left( \tilde{u}\mid h\left( g\cdot \right) \right) =\left( 
\tilde{u}\left( g^{-1}\cdot \right) \mid h\right) $ and $T\left( u\right)
h\left( g\cdot \right) =T\left( u\right) h$, so that $\left( \tilde{u}\left(
g^{-1}\cdot \right) \mid h\right) =\left( \tilde{u}\mid h\right) $. This
means $\tilde{u}\left( g^{-1}\cdot \right) =\tilde{u}$ for all $g\in O\left(
N\right) $, i.e., $\tilde{u}\in H_{V,\mathrm{r}}^{1}$. Now assume $I^{\prime
}\left( u\right) =0$ in the dual space of $H_{V,\mathrm{r}}^{1}$. Then we
have $\left( \tilde{u}\mid h\right) =T\left( u\right) h=I^{\prime }\left(
u\right) h=0$ for all $h\in H_{V,\mathrm{r}}^{1}$, which implies $\tilde{u}%
=0 $. This gives $T\left( u\right) h=0$ for all $h\in H_{V}^{1}$, which is
the thesis (\ref{PROP:symm-crit: th}).%
\endproof
\bigskip

For future reference, we point out here that, by assumption $\left(
h_{1}\right) $, if $\left( f_{q_{1},q_{2}}\right) $ holds then $\exists 
\tilde{M}>0$ such that for almost every $r>0$ and all $t\in \mathbb{R}$ one has 
\begin{equation}
\left| G\left( r,t\right) -g\left( r,0\right) t\right| \leq \tilde{M}K\left(
r\right) \min \left\{ \left| t\right| ^{q_{1}},\left| t\right|
^{q_{2}}\right\} .  \label{G_pq}
\end{equation}

\begin{lemma}
\label{LEM:MPgeom}Let $L_{0}$ be the norm of the operator of assumption $%
\left( h_{3}\right) $. If there exist $q_{1},q_{2}>1$ such that $\left(
f_{q_{1},q_{2}}\right) $ and $\left( \mathcal{S}_{q_{1},q_{2}}^{\prime
}\right) $ hold, then there exist two constants $c_{1},c_{2}>0$ such that 
\begin{equation}
I\left( u\right) \geq \frac{1}{2}\left\| u\right\| ^{2}-c_{1}\left\|
u\right\| ^{q_{1}}-c_{2}\left\| u\right\| ^{q_{2}}-L_{0}\left\| u\right\|
\qquad \text{for all }u\in H_{V,\mathrm{r}}^{1}.  \label{LEM:MPgeom: th}
\end{equation}
If $\left( \mathcal{S}_{q_{1},q_{2}}^{\prime \prime }\right) $ also holds,
then $\forall \varepsilon >0$ there exist two constants $c_{1}\left(
\varepsilon \right) ,c_{2}\left( \varepsilon \right) >0$ such that (\ref
{LEM:MPgeom: th}) holds both with $c_{1}=\varepsilon $, $c_{2}=c_{2}\left(
\varepsilon \right) $ and with $c_{1}=c_{1}\left( \varepsilon \right) $, $%
c_{2}=\varepsilon $.
\end{lemma}

\proof 
Let $i\in \left\{ 1,2\right\} $. By the monotonicity of $\mathcal{S}_{0}$
and $\mathcal{S}_{\infty }$, it is not restrictive to assume $R_{1}<R_{2}$
in hypothesis $\left( \mathcal{S}_{q_{1},q_{2}}^{\prime }\right) $. Then, by 
\cite[Lemma 1]{BGR-p1} and the continuous embedding $H_{V}^{1}%
\hookrightarrow L_{\mathrm{loc}}^{2}(\mathbb{R}^{N})$, there exists a constant $%
c_{R_{1},R_{2}}^{\left( i\right) }>0$ such that for all $u\in H_{V,\mathrm{r}%
}^{1}$ we have 
\[
\int_{B_{R_{2}}\setminus B_{R_{1}}}K\left( \left| x\right| \right) \left|
u\right| ^{q_{i}}dx\leq c_{R_{1},R_{2}}^{\left( i\right) }\left\| u\right\|
^{q_{i}}. 
\]
Therefore, by (\ref{G_pq}) and the definitions of $\mathcal{S}_{0}$ and $%
\mathcal{S}_{\infty }$, we obtain 
\begin{eqnarray}
&&\left| \int_{\mathbb{R}^{N}}G\left( \left| x\right| ,u\right) dx\right| 
\nonumber \\
&\leq &\int_{\mathbb{R}^{N}}\left| G\left( \left| x\right| ,u\right) -g\left(
\left| x\right| ,0\right) u\right| dx+\left| \int_{\mathbb{R}^{N}}g\left(
\left| x\right| ,0\right) u\,dx\right|  \nonumber \\
&\leq &\tilde{M}\int_{\mathbb{R}^{N}}K\left( \left| x\right| \right) \min
\left\{ \left| u\right| ^{q_{1}},\left| u\right| ^{q_{2}}\right\}
dx+L_{0}\left\| u\right\|  \nonumber \\
&\leq &\tilde{M}\left( \int_{B_{R_{1}}}K\left( \left| x\right| \right)
\left| u\right| ^{q_{1}}dx+\int_{B_{R_{2}}^{c}}K\left( \left| x\right|
\right) \left| u\right| ^{q_{2}}dx+\int_{B_{R_{2}}\setminus
B_{R_{1}}}K\left( \left| x\right| \right) \left| u\right| ^{q_{i}}dx\right)
+L_{0}\left\| u\right\|  \nonumber \\
&\leq &\tilde{M}\left( \left\| u\right\| ^{q_{1}}\mathcal{S}_{0}\left(
q_{1},R_{1}\right) +\left\| u\right\| ^{q_{2}}\mathcal{S}_{\infty }\left(
q_{2},R_{2}\right) +c_{R_{1},R_{2}}^{\left( i\right) }\left\| u\right\|
^{q_{i}}\right) +L_{0}\left\| u\right\|  \label{LEM:MPgeom: pf} \\
&=&c_{1}\left\| u\right\| ^{q_{1}}+c_{2}\left\| u\right\|
^{q_{2}}+L_{0}\left\| u\right\| ,  \nonumber
\end{eqnarray}
with obvious definition of the constants $c_{1}$ and $c_{2}$, independent of 
$u$. This yields (\ref{LEM:MPgeom: th}). If $\left( \mathcal{S}%
_{q_{1},q_{2}}^{\prime \prime }\right) $ also holds, then $\forall
\varepsilon >0$ we can fix $R_{1,\varepsilon }<R_{2,\varepsilon }$ such that 
$\tilde{M}\mathcal{S}_{0}\left( q_{1},R_{1,\varepsilon }\right) <\varepsilon 
$ and $\tilde{M}\mathcal{S}_{\infty }\left( q_{2},R_{2,\varepsilon }\right)
<\varepsilon $, so that inequality (\ref{LEM:MPgeom: pf}) becomes 
\[
\left| \int_{\mathbb{R}^{N}}G\left( \left| x\right| ,u\right) dx\right| \leq
\varepsilon \left\| u\right\| ^{q_{1}}+\varepsilon \left\| u\right\|
^{q_{2}}+c_{R_{1,\varepsilon },R_{2,\varepsilon }}^{\left( i\right) }\left\|
u\right\| ^{q_{i}}+L_{0}\left\| u\right\| . 
\]
The result then ensues by taking $i=2$ and $c_{2}\left( \varepsilon \right)
=\varepsilon +c_{R_{1,\varepsilon },R_{2,\varepsilon }}^{\left( 2\right) }$,
or $i=1$ and $c_{1}\left( \varepsilon \right) =\varepsilon
+c_{R_{1,\varepsilon },R_{2,\varepsilon }}^{\left( 1\right) }$.%
\endproof
\bigskip

Henceforth, we will assume that the hypotheses of Theorem \ref{THM:ex} also
include the following condition: 
\begin{equation}
g\left( r,t\right) =0\quad \text{for all }r>0~\text{and~}t<0.  \label{g=0}
\end{equation}
This can be done without restriction, since the theorem concerns \emph{%
nonnegative} critical points and all its assumptions still hold true if we
replace $g\left( r,t\right) $ with $g\left( r,t\right) \chi _{\mathbb{R}%
_{+}}\left( t\right) $ ($\chi _{\mathbb{R}_{+}}$ is the characteristic function
of $\mathbb{R}_{+}$).

\begin{lemma}
\label{LEM:PS}Under the assumptions of each of Theorems \ref{THM:ex}
(including (\ref{g=0})) and \ref{THM:mult}, the functional $I:H_{V,\mathrm{r}%
}^{1}\rightarrow \mathbb{R}$ satisfies the Palais-Smale condition.
\end{lemma}

\proof 
By (\ref{g=0}) and $\left( g_{5}\right) $ respectively, under the
assumptions of each of Theorems \ref{THM:ex} and \ref{THM:mult} we have that
either $g$ satisfies $\left( g_{1}\right) $ for all $t\in \mathbb{R}$, or $%
K\left( \left| \cdot \right| \right) \in L^{1}(\mathbb{R}^{N})$ and $g$
satisfies 
\begin{equation}
\theta G\left( r,t\right) \leq g\left( r,t\right) t\quad \text{for almost
every }r>0~\text{and all~}\left| t\right| \geq t_{0}.  \label{LEM:PS: AR}
\end{equation}
Let $\left\{ u_{n}\right\} $ be a sequence in $H_{V,\mathrm{r}}^{1}$ such
that $\left\{ I\left( u_{n}\right) \right\} $ is bounded and $I^{\prime
}\left( u_{n}\right) \rightarrow 0$ in the dual space of $H_{V,\mathrm{r}%
}^{1}$. Hence 
\[
\frac{1}{2}\left\| u_{n}\right\| ^{2}-\int_{\mathbb{R}^{N}}G\left( \left|
x\right| ,u_{n}\right) dx=O\left( 1\right) \quad \text{and}\quad \left\|
u_{n}\right\| ^{2}-\int_{\mathbb{R}^{N}}g\left( \left| x\right| ,u_{n}\right)
u_{n}dx=o\left( 1\right) \left\| u_{n}\right\| . 
\]
If $g$ satisfies $\left( g_{1}\right) $, then we get 
\[
\frac{1}{2}\left\| u_{n}\right\| ^{2}+O\left( 1\right) =\int_{\mathbb{R}%
^{N}}G\left( \left| x\right| ,u_{n}\right) dx\leq \frac{1}{\theta }\int_{%
\mathbb{R}^{N}}g\left( \left| x\right| ,u_{n}\right) u_{n}dx=\frac{1}{\theta }%
\left\| u_{n}\right\| ^{2}+o\left( 1\right) \left\| u_{n}\right\| , 
\]
which implies that $\left\{ \left\| u_{n}\right\| \right\} $ is bounded
since $\theta >2$. If $K\left( \left| \cdot \right| \right) \in L^{1}(\mathbb{R}%
^{N})$ and $g$ satisfies (\ref{LEM:PS: AR}), then we slightly modify the
argument: we have 
\begin{eqnarray*}
\int_{\left\{ \left| u_{n}\right| \geq t_{0}\right\} }g\left( \left|
x\right| ,u_{n}\right) u_{n}dx &=&\int_{\mathbb{R}^{N}}g\left( \left| x\right|
,u_{n}\right) u_{n}dx-\int_{\left\{ \left| u_{n}\right| <t_{0}\right\}
}g\left( \left| x\right| ,u_{n}\right) u_{n}dx \\
&\leq &\int_{\mathbb{R}^{N}}g\left( \left| x\right| ,u_{n}\right)
u_{n}dx+\int_{\left\{ \left| u_{n}\right| <t_{0}\right\} }\left| g\left(
\left| x\right| ,u_{n}\right) u_{n}\right| dx
\end{eqnarray*}
where (thanks to $\left( h_{1}\right) $ and $\left( f_{q_{1},q_{2}}\right) $%
) 
\begin{eqnarray*}
\int_{\left\{ \left| u_{n}\right| <t_{0}\right\} }\left| g\left( \left|
x\right| ,u_{n}\right) u_{n}\right| dx &\leq &\int_{\left\{ \left|
u_{n}\right| <t_{0}\right\} }K\left( \left| x\right| \right) \left| f\left(
u_{n}\right) \right| \left| u_{n}\right| dx \\
&\leq &M\int_{\left\{ \left| u_{n}\right| <t_{0}\right\} }K\left( \left|
x\right| \right) \min \left\{ \left| u_{n}\right| ^{q_{1}},\left|
u_{n}\right| ^{q_{2}}\right\} dx \\
&\leq &M\min \left\{ t_{0}^{q_{1}},t_{0}^{q_{2}}\right\} \int_{\left\{
\left| u_{n}\right| <t_{0}\right\} }K\left( \left| x\right| \right) dx\leq
M\min \left\{ t_{0}^{q_{1}},t_{0}^{q_{2}}\right\} \left\| K\right\| _{L^{1}(%
\mathbb{R}^{N})},
\end{eqnarray*}
so that, by (\ref{G_pq}), we obtain 
\begin{eqnarray*}
\frac{1}{2}\left\| u_{n}\right\| ^{2}+O\left( 1\right) &=&\int_{\mathbb{R}%
^{N}}G\left( \left| x\right| ,u_{n}\right) dx=\int_{\left\{ \left|
u_{n}\right| <t_{0}\right\} }G\left( \left| x\right| ,u_{n}\right)
dx+\int_{\left\{ \left| u_{n}\right| \geq t_{0}\right\} }G\left( \left|
x\right| ,u_{n}\right) dx \\
&\leq &\tilde{M}\int_{\left\{ \left| u_{n}\right| <t_{0}\right\} }K\left(
\left| x\right| \right) \min \left\{ \left| u_{n}\right| ^{q_{1}},\left|
u_{n}\right| ^{q_{2}}\right\} dx+\frac{1}{\theta }\int_{\left\{ \left|
u_{n}\right| \geq t_{0}\right\} }g\left( \left| x\right| ,u_{n}\right)
u_{n}dx \\
&\leq &\tilde{M}\min \left\{ t_{0}^{q_{1}},t_{0}^{q_{2}}\right\} \left\|
K\right\| _{L^{1}(\mathbb{R}^{N})}+\frac{1}{\theta }\int_{\mathbb{R}^{N}}g\left(
\left| x\right| ,u_{n}\right) u_{n}dx+\frac{M}{\theta }\min \left\{
t_{0}^{q_{1}},t_{0}^{q_{2}}\right\} \left\| K\right\| _{L^{1}(\mathbb{R}^{N})}
\\
&=&\left( \tilde{M}+\frac{M}{\theta }\right) \min \left\{
t_{0}^{q_{1}},t_{0}^{q_{2}}\right\} \left\| K\right\| _{L^{1}(\mathbb{R}^{N})}+%
\frac{1}{\theta }\left\| u_{n}\right\| ^{2}+o\left( 1\right) \left\|
u_{n}\right\| .
\end{eqnarray*}
This yields again that $\left\{ \left\| u_{n}\right\| \right\} $ is bounded.
Now, since the embedding $H_{V,\mathrm{r}}^{1}\hookrightarrow
L_{K}^{q_{1}}+L_{K}^{q_{2}}$ is compact by assumption $\left( \mathcal{S}%
_{q_{1},q_{2}}^{\prime \prime }\right) $ (see \cite[Theorem 1]{BGR-p1}) and
the operator $u\mapsto \int_{\mathbb{R}^{N}}G\left( \left| x\right| ,u\right)
dx $ is of class $C^{1}$ on $L_{K}^{q_{1}}+L_{K}^{q_{2}}$ (see the proof of
Proposition \ref{PROP:diff} above), it is a standard exercise to conclude
that $\left\{ u_{n}\right\} $ has a strongly convergent subsequence in $H_{V,%
\mathrm{r}}^{1}$.%
\endproof
\bigskip

\noindent \textbf{Proof of Theorem \ref{THM:ex}.}\quad 
We want to apply the Mountain-Pass Theorem \cite{Ambr-Rab}. To this end, from (\ref{LEM:MPgeom:
th}) of Lemma \ref{LEM:MPgeom} we deduce that, since $L_{0}=0$ and $%
q_{1},q_{2}>2$, there exists $\rho >0$ such that 
\begin{equation}
\inf_{u\in H_{V,\mathrm{r}}^{1},\,\left\| u\right\| =\rho }I\left( u\right)
>0=I\left( 0\right) .  \label{mp-geom}
\end{equation}
Therefore, taking into account Proposition \ref{PROP:diff} and Lemma \ref
{LEM:PS}, we need only to check that $\exists \bar{u}\in H_{V,\mathrm{r}%
}^{1} $ such that $\left\| \bar{u}\right\| >\rho $ and $I\left( \bar{u}%
\right) <0$. In order to check this, from assumption $\left( g_{3}\right) $
(which holds in any case, according to Remark \ref{RMK:thm:ex}.\ref
{RMK:thm:ex-1}), we infer that 
\[
G\left( r,t\right) \geq \frac{G\left( r,t_{0}\right) }{t_{0}^{\theta }}%
t^{\theta }\text{\quad for almost every }r>0\text{~and all }t\geq t_{0}. 
\]
Then, by assumption $\left( h_{0}\right) $, we fix a nonnegative function $%
u_{0}\in C_{c}^{\infty }(B_{r_{2}}\setminus \overline{B}_{r_{1}})\cap H_{V,%
\mathrm{r}}^{1}$ such that the set $\{x\in \mathbb{R}^{N}:u_{0}\left( x\right)
\geq t_{0}\}$ has positive Lebesgue measure. We now distinguish the case of
assumptions $\left( g_{1}\right) $ and $\left( g_{2}\right) $ from the case
of $K\left( \left| \cdot \right| \right) \in L^{1}(\mathbb{R}^{N})$. In the
first one, $\left( g_{1}\right) $ and $\left( g_{2}\right) $ ensure that $%
G\geq 0$ and $G\left( \cdot ,t_{0}\right) >0$ almost everywhere, so that for
every $\lambda >1$ we get 
\begin{eqnarray*}
\int_{\mathbb{R}^{N}}G\left( \left| x\right| ,\lambda u_{0}\right) dx &\geq
&\int_{\left\{ \lambda u_{0}\geq t_{0}\right\} }G\left( \left| x\right|
,\lambda u_{0}\right) dx\geq \frac{\lambda ^{\theta }}{t_{0}^{\theta }}%
\int_{\left\{ \lambda u_{0}\geq t_{0}\right\} }G\left( \left| x\right|
,t_{0}\right) u_{0}^{\theta }dx \\
&\geq &\frac{\lambda ^{\theta }}{t_{0}^{\theta }}\int_{\left\{ u_{0}\geq
t_{0}\right\} }G\left( \left| x\right| ,t_{0}\right) u_{0}^{\theta }dx\geq
\lambda ^{\theta }\int_{\left\{ u_{0}\geq t_{0}\right\} }G\left( \left|
x\right| ,t_{0}\right) dx>0.
\end{eqnarray*}
Since $\theta >2$, this gives 
\[
\lim_{\lambda \rightarrow +\infty }I\left( \lambda u_{0}\right) \leq
\lim_{\lambda \rightarrow +\infty }\left( \frac{\lambda ^{2}}{2}\left\|
u_{0}\right\| ^{2}-\lambda ^{\theta }\int_{\left\{ u_{0}\geq t_{0}\right\}
}G\left( \left| x\right| ,t_{0}\right) dx\right) =-\infty . 
\]
If $K\left( \left| \cdot \right| \right) \in L^{1}(\mathbb{R}^{N})$, assumption 
$\left( g_{3}\right) $ still gives $G\left( \cdot ,t_{0}\right) >0$ almost
everywhere and from (\ref{G_pq}) we infer that 
\[
G\left( r,t\right) \geq -\tilde{M}K\left( r\right) \min \left\{
t_{0}^{q_{1}},t_{0}^{q_{2}}\right\} \text{\quad for almost every }r>0\text{%
~and all }0\leq t\leq t_{0}. 
\]
Therefore, arguing as before about the integral over $\left\{ \lambda
u_{0}\geq t_{0}\right\} $, for every $\lambda >1$ we obtain 
\begin{eqnarray*}
\int_{\mathbb{R}^{N}}G\left( \left| x\right| ,\lambda u_{0}\right) dx
&=&\int_{\left\{ \lambda u_{0}<t_{0}\right\} }G\left( \left| x\right|
,\lambda u_{0}\right) dx+\int_{\left\{ \lambda u_{0}\geq t_{0}\right\}
}G\left( \left| x\right| ,\lambda u_{0}\right) dx \\
&\geq &-\tilde{M}\min \left\{ t_{0}^{q_{1}},t_{0}^{q_{2}}\right\}
\int_{\left\{ \lambda u_{0}<t_{0}\right\} }K\left( \left| x\right| \right)
dx+\lambda ^{\theta }\int_{\left\{ u_{0}\geq t_{0}\right\} }G\left( \left|
x\right| ,t_{0}\right) dx,
\end{eqnarray*}
which implies 
\[
\lim_{\lambda \rightarrow +\infty }I\left( \lambda u_{0}\right) \leq
\lim_{\lambda \rightarrow +\infty }\left( \frac{\lambda ^{2}}{2}\left\|
u_{0}\right\| ^{2}+\tilde{M}\min \left\{ t_{0}^{q_{1}},t_{0}^{q_{2}}\right\}
\left\| K\right\| _{L^{1}(\mathbb{R}^{N})}-\lambda ^{\theta }\int_{\left\{
u_{0}\geq t_{0}\right\} }G\left( \left| x\right| ,t_{0}\right) dx\right)
=-\infty . 
\]
So, in any case, we can take $\bar{u}=\lambda u_{0}$ with $\lambda $
sufficiently large and the Mountain-Pass Theorem provides the existence of a
nonzero critical point $u\in H_{V,\mathrm{r}}^{1}$ for $I$. Since (\ref{g=0}%
) implies $I^{\prime }\left( u\right) u_{-}=-\left\| u_{-}\right\| ^{2}$
(where $u_{-}\in H_{V,\mathrm{r}}^{1}$ is the negative part of $u$), one
concludes that $u_{-}=0$, i.e., $u$ is nonnegative.%
\endproof
\bigskip

\noindent \textbf{Proof of Theorem \ref{THM:mult}.}\quad 
By the oddness assumption $\left( g_{5}\right) $, one has $I\left( u\right) =I\left(
-u\right) $ for all $u\in H_{V,\mathrm{r}}^{1}$ and thus we can apply the
Symmetric Mountain-Pass Theorem (see e.g. \cite[Chapter 1]{Rabi}). To this
end, we deduce (\ref{mp-geom}) as in the proof of Theorem \ref{THM:ex} and
therefore, thanks to Proposition \ref{PROP:diff} and Lemma \ref{LEM:PS}, we
need only to show that $I$ satisfies the following geometrical condition:
for any finite dimensional subspace $Y\neq \left\{ 0\right\} $ of $H_{V,%
\mathrm{r}}^{1}$ there exists $R>0$ such that $I\left( u\right) \leq 0$ for
all $u\in Y$ with $\left\| u\right\| \geq R$. In fact, it is sufficient to
prove that any diverging sequence in $Y$ admits a subsequence on which $I$
is nonpositive. So, let $\left\{ u_{n}\right\} \subseteq Y$ be such that $%
\left\| u_{n}\right\| \rightarrow +\infty $. Since all norms are equivalent
on $Y$, by (\ref{Vu}) one has 
\begin{equation}
\left\| u_{n}\right\| _{L_{K}^{p}(\Lambda _{u_{n}})}+\left\| u_{n}\right\|
_{L_{K}^{q}(\Lambda _{u_{n}}^{c})}\geq \left\| u_{n}\right\|
_{L_{K}^{q_{1}}+L_{K}^{q_{2}}}\geq m_{1}\left\| u_{n}\right\| \rightarrow
+\infty  \label{THM2-pf: eq_norms}
\end{equation}
for some constant $m_{1}>0$, where $p:=\min \left\{ q_{1},q_{2}\right\} $
and $q:=\max \left\{ q_{1},q_{2}\right\} $. Hence, up to a subsequence, at
least one of the sequences $\{\left\| u_{n}\right\| _{L_{K}^{p}(\Lambda
_{u_{n}})}\}$, $\{\left\| u_{n}\right\| _{L_{K}^{q}(\Lambda _{u_{n}}^{c})}\}$
diverges. We now use assumptions $\left( g_{4}\right) $ and $\left(
g_{5}\right) $ to deduce that 
\[
G\left( r,t\right) \geq mK\left( r\right) \min \left\{ \left| t\right|
^{q_{1}},\left| t\right| ^{q_{2}}\right\} \quad \text{for almost every }r>0~%
\text{and all }t\in \mathbb{R}, 
\]
which implies 
\begin{eqnarray*}
\int_{\mathbb{R}^{N}}G\left( \left| x\right| ,u_{n}\right) dx &\geq &m\int_{%
\mathbb{R}^{N}}K\left( \left| x\right| \right) \min \left\{ \left| u_{n}\right|
^{q_{1}},\left| u_{n}\right| ^{q_{2}}\right\} dx \\
&=&m\int_{\Lambda _{u_{n}}}K\left( \left| x\right| \right) \left|
u_{n}\right| ^{p}dx+m\int_{\Lambda _{u_{n}}^{c}}K\left( \left| x\right|
\right) \left| u_{n}\right| ^{q}dx.
\end{eqnarray*}
Hence, using inequalities (\ref{THM2-pf: eq_norms}), there exists a constant 
$m_{2}>0$ such that 
\[
I\left( u_{n}\right) \leq m_{2}\left( \left\| u_{n}\right\|
_{L_{K}^{p}(\Lambda _{u_{n}})}^{2}+\left\| u_{n}\right\| _{L_{K}^{q}(\Lambda
_{u_{n}}^{c})}^{2}\right) -m\left( \left\| u_{n}\right\| _{L_{K}^{p}(\Lambda
_{u_{n}})}^{p}+\left\| u_{n}\right\| _{L_{K}^{q}(\Lambda
_{u_{n}}^{c})}^{q}\right) , 
\]
so that $I\left( u_{n}\right) \rightarrow -\infty $ since $p,q>2$. The
Symmetric Mountain-Pass Theorem thus implies the existence of an unbounded
sequence of critical values for $I$ and this completes the proof.%
\endproof
\bigskip

\begin{lemma}
\label{LEM:bdd}Under the assumptions of each of Theorems \ref{THM:ex-sub}
and \ref{THM:mult-sub}, the functional $I:H_{V,\mathrm{r}}^{1}\rightarrow 
\mathbb{R}$ is bounded from below and coercive. In particular, if $g$ satisfies 
$\left( g_{6}\right) $, then 
\begin{equation}
\inf_{v\in H_{V,\mathrm{r}}^{1}}I\left( v\right) <0.  \label{LEM:bdd: inf<0}
\end{equation}
\end{lemma}

\proof 
The fact that $I$ is bounded below and coercive on $H_{V,\mathrm{r}}^{1}$ is
a consequence of Lemma \ref{LEM:MPgeom}. Indeed, the result readily follows
from (\ref{LEM:MPgeom: th}) if $q_{1},q_{2}\in \left( 1,2\right) $, while,
if $\max \left\{ q_{1},q_{2}\right\} =2>\min \left\{ q_{1},q_{2}\right\} >1$%
, we fix $\varepsilon <1/2$ and use the second part of the lemma in order to
get 
\[
I\left( u\right) \geq \left( \frac{1}{2}-\varepsilon \right) \left\|
u\right\| ^{2}-c\left( \varepsilon \right) \left\| u\right\| ^{\min \left\{
q_{1},q_{2}\right\} }-L_{0}\left\| u\right\| \qquad \text{for all }u\in H_{V,%
\mathrm{r}}^{1}, 
\]
which yields again the conclusion. In order to prove (\ref{LEM:bdd: inf<0}),
we use assumption $\left( h_{0}\right) $ to fix a function $u_{0}\in
C_{c}^{\infty }(B_{r_{2}}\setminus \overline{B}_{r_{1}})\cap H_{V,\mathrm{r}%
}^{1}$ such that $0\leq u_{0}\leq t_{0}$, $u_{0}\neq 0$. Then, by assumption 
$\left( g_{6}\right) $, for every $0<\lambda <1$ we get that $\lambda
u_{0}\in H_{V,\mathrm{r}}^{1}$ satisfies 
\[
I\left( \lambda u_{0}\right) =\frac{1}{2}\left\| \lambda u_{0}\right\|
^{2}-\int_{\mathbb{R}^{N}}G\left( \left| x\right| ,\lambda u_{0}\right) dx\leq 
\frac{\lambda ^{2}}{2}\left\| u_{0}\right\| ^{2}-\lambda ^{\theta }m\int_{%
\mathbb{R}^{N}}K\left( \left| x\right| \right) u_{0}^{\theta }dx. 
\]
Since $\theta <2$, this implies $I\left( \lambda u_{0}\right) <0$ for $%
\lambda $ sufficiently small and therefore (\ref{LEM:bdd: inf<0}) ensues.%
\endproof
\bigskip

\noindent \textbf{Proof of Theorem \ref{THM:ex-sub}.}\quad 
Let 
\[
\mu :=\inf_{v\in H_{V,\mathrm{r}}^{1}}I\left( v\right) 
\]
and take any minimizing sequence $\left\{ v_{n}\right\} $ for $\mu $. From
Lemma \ref{LEM:bdd} we have that the functional $I:H_{V,\mathrm{r}%
}^{1}\rightarrow \mathbb{R}$ is bounded from below and coercive, so that $\mu
\in \mathbb{R}$ and $\left\{ v_{n}\right\} $ is bounded in $H_{V,\mathrm{r}%
}^{1} $. Thanks to assumption $\left( \mathcal{S}_{q_{1},q_{2}}^{\prime
\prime }\right) $, the embedding $H_{V,\mathrm{r}}^{1}\hookrightarrow
L_{K}^{q_{1}}+L_{K}^{q_{2}}$ is compact (see \cite[Theorem 1]{BGR-p1}) and
thus we can assume that there exists $u\in H_{V,\mathrm{r}}^{1}$ such that,
up to a subsequence, one has: 
\[
\begin{array}{ll}
v_{n}\rightharpoonup u & \text{in }H_{V,\mathrm{r}}^{1},\medskip \\ 
v_{n}\rightarrow u & \text{in }L_{K}^{q_{1}}+L_{K}^{q_{2}}.
\end{array}
\]
Then, thanks to $\left( h_{3}\right) $ and the continuity of the functional (%
\ref{NemOp}) on $L_{K}^{q_{1}}+L_{K}^{q_{2}}$ (see the proof of Proposition 
\ref{PROP:diff} above), $u$ satisfies 
\[
\int_{\mathbb{R}^{N}}G\left( \left| x\right| ,v_{n}\right) dx=\int_{\mathbb{R}%
^{N}}\left( G\left( \left| x\right| ,v_{n}\right) -g\left( \left| x\right|
,0\right) v_{n}\right) dx+\int_{\mathbb{R}^{N}}g\left( \left| x\right|
,0\right) v_{n}dx\rightarrow \int_{\mathbb{R}^{N}}G\left( \left| x\right|
,u\right) dx. 
\]
By the weak lower semi-continuity of the norm, this implies 
\[
I\left( u\right) =\frac{1}{2}\left\| u\right\| ^{2}-\int_{\mathbb{R}%
^{N}}G\left( \left| x\right| ,u\right) dx\leq \lim_{n\rightarrow \infty
}\left( \frac{1}{2}\left\| v_{n}\right\| ^{2}-\int_{\mathbb{R}^{N}}G\left(
\left| x\right| ,v_{n}\right) dx\right) =\mu 
\]
and thus we conclude $I\left( u\right) =\mu $. It remains to show that $%
u\neq 0$. If $g$ satisfies $\left( g_{6}\right) $, then we have $\mu <0$ by
Lemma \ref{LEM:bdd} and therefore it must be $u\neq 0$, since $I\left(
0\right) =0$. If $\left( g_{7}\right) $ holds, assume by contradiction that $%
u=0$. Since $u$ is a critical point of $I\in C^{1}(H_{V,\mathrm{r}}^{1};\mathbb{%
R})$, from (\ref{PROP:diff: I'(u)h=}) we get 
\[
\int_{\mathbb{R}^{N}}g\left( \left| x\right| ,0\right) h\,dx=0,\quad \forall
h\in C_{c,\mathrm{rad}}^{\infty }(B_{r_{2}}\setminus \overline{B}%
_{r_{1}})\subset H_{V,\mathrm{r}}^{1}. 
\]
This implies $g\left( \cdot ,0\right) =0$ almost everywhere in $\left(
r_{1},r_{2}\right) $, which is a contradiction.%
\endproof
\bigskip

\noindent \textbf{Proof of Corollary \ref{COR:ex-sub}.}\quad 
Setting 
\[
\widetilde{g}\left( r,t\right) :=\left\{ 
\begin{array}{lll}
g\left( r,t\right) &  & \text{if }t\geq 0 \\ 
2g\left( r,0\right) -g\left( r,\left| t\right| \right) &  & \text{if }t<0
\end{array}
\right. 
\]
and 
\[
\widetilde{G}\left( r,t\right) :=\int_{0}^{t}\widetilde{g}\left( r,s\right)
ds=\left\{ 
\begin{array}{lll}
G\left( r,t\right) &  & \text{if }t\geq 0 \\ 
2g\left( r,0\right) t+G\left( r,\left| t\right| \right) &  & \text{if }t<0,
\end{array}
\right. 
\]
it is easy to check that the function $\widetilde{g}$ still satisfies all
the assumptions of Theorem \ref{THM:ex-sub}. We just observe that $%
\widetilde{g}$ satisfies $\left( g_{6}\right) $ or $\left( g_{7}\right) $ if
so does $g$, and that for almost every $r>0$ and all $t\in \mathbb{R}$ one has 
\[
\left| \widetilde{g}\left( r,t\right) -\widetilde{g}\left( r,0\right)
\right| =\left| g\left( r,\left| t\right| \right) -g\left( r,0\right)
\right| \leq K\left( r\right) \left| f\left( \left| t\right| \right) \right|
\quad \text{with}\quad \left| f\left( \left| t\right| \right) \right| \leq
M\min \left\{ \left| t\right| ^{q_{1}-1},\left| t\right| ^{q_{2}-1}\right\}
. 
\]
Then, by Theorem \ref{THM:ex-sub}, there exists $\widetilde{u}\neq 0$ such
that 
\[
\widetilde{I}\left( \widetilde{u}\right) =\min_{u\in H_{V,\mathrm{r}}^{1}}%
\widetilde{I}\left( u\right) ,\quad \text{where}\quad \widetilde{I}\left(
u\right) :=\frac{1}{2}\left\| u\right\| ^{2}-\int_{\mathbb{R}^{N}}\widetilde{G}%
\left( \left| x\right| ,u\right) dx. 
\]
For every $u\in H_{V,\mathrm{r}}^{1}$ one has 
\begin{eqnarray}
\widetilde{I}\left( u\right) &=&\frac{1}{2}\left\| u\right\|
^{2}-\int_{\left\{ u\geq 0\right\} }G\left( \left| x\right| ,u\right)
dx-2\int_{\left\{ u<0\right\} }g\left( \left| x\right| ,0\right)
u\,dx-\int_{\left\{ u<0\right\} }G\left( \left| x\right| ,\left| u\right|
\right) dx  \nonumber \\
&=&\frac{1}{2}\left\| u\right\| ^{2}-\int_{\mathbb{R}^{N}}G\left( \left|
x\right| ,\left| u\right| \right) dx+2\int_{\mathbb{R}^{N}}g\left( \left|
x\right| ,0\right) u_{-}\,dx  \label{I-tilde=} \\
&=&I\left( \left| u\right| \right) +2\int_{\mathbb{R}^{N}}g\left( \left|
x\right| ,0\right) u_{-}\,dx,  \nonumber
\end{eqnarray}
which implies that $\widetilde{u}$ satisfies (\ref{COR:ex-sub: th}), as one
readily checks that 
\[
\inf_{u\in H_{V,\mathrm{r}}^{1}}\left( I\left( \left| u\right| \right)
+2\int_{\mathbb{R}^{N}}g\left( \left| x\right| ,0\right) u_{-}\,dx\right)
=\inf_{u\in H_{V,\mathrm{r}}^{1},\,u\geq 0}I\left( u\right) . 
\]
Moreover, since $G\left( r,\left| t\right| \right) =\widetilde{G}\left(
r,\left| t\right| \right) $ and $g\left( \cdot ,0\right) \geq 0$, (\ref
{I-tilde=}) gives 
\[
\widetilde{I}\left( u\right) =\widetilde{I}\left( \left| u\right| \right)
+2\int_{\mathbb{R}^{N}}g\left( \left| x\right| ,0\right) u_{-}\,dx\geq 
\widetilde{I}\left( \left| u\right| \right) 
\]
and hence $\left| \widetilde{u}\right| \in H_{V,\mathrm{r}}^{1}$ is still a
minimizer for $\widetilde{I}$, so that we can assume $\widetilde{u}\geq 0$.
Finally, $\widetilde{u}$ is a critical point for $I$ since $\widetilde{u}$
is a critical point of $\widetilde{I}$ and $\widetilde{g}\left( r,t\right)
=g\left( r,t\right) $ for avery $t\geq 0$. 
\endproof
\bigskip

In proving Theorem \ref{THM:mult-sub}, we will use a well known abstract
result from \cite{Clark,Heinz}. We recall it here in a version given in \cite
{Wang01}.

\begin{theorem}[{\cite[Lemma 2.4]{Wang01}}]
\label{THM:wang}
Let $X$ be a real Banach space and let $J\in C^{1}(X;\mathbb{R})$. Assume that $J$ satisfies the Palais-Smale
condition, is even, bounded from below and such that $J\left( 0\right) =0$.
Assume furthermore that $\forall k\in \mathbb{N}\setminus \left\{ 0\right\} $
there exist $\rho _{k}>0$ and a $k$-dimensional subspace $X_{k}$ of $X$ such
that 
\begin{equation}
\sup_{u\in X_{k},\,\left\| u\right\| _{X}=\rho _{k}}J\left( u\right) <0.
\label{THM:wang: geom}
\end{equation}
Then $J$ has a sequence of critical values $c_{k}<0$ such that $%
\displaystyle%
\lim_{k\rightarrow \infty }c_{k}=0$.
\end{theorem}

\noindent \textbf{Proof of Theorem \ref{THM:mult-sub}.}\quad 
Since $I:H_{V,\mathrm{r}}^{1}\rightarrow \mathbb{R}$ satisfies $I\left( 0\right) =0$ 
and is of class $C^{1}$ by Proposition \ref{PROP:diff}, even by assumption $\left(
g_{5}\right) $ and bounded below by Lemma \ref{LEM:bdd}, for applying
Theorem \ref{THM:wang} (with $X=H_{V,\mathrm{r}}^{1}$ and $J=I$) we need
only to show that $I$ satisfies the Palais-Smale condition and the geometric
condition (\ref{THM:wang: geom}). By coercivity (Lemma \ref{LEM:bdd}), every
Palais-Smale sequence for $I$ is bounded in $H_{V,\mathrm{r}}^{1}$ and one
obtains the existence of a strongly convergent subsequence as in the proof
of Lemma \ref{LEM:PS}. In order to check (\ref{THM:wang: geom}), we first
deduce from $\left( g_{5}\right) $ and $\left( g_{6}\right) $ that 
\begin{equation}
G\left( r,t\right) \geq mK\left( r\right) \left| t\right| ^{\theta }\quad 
\text{for almost every }r>0~\text{and all }\left| t\right| \leq t_{0}.
\label{G>}
\end{equation}
Then, for any $k\in \mathbb{N}\setminus \left\{ 0\right\} $, we take $k$
linearly independent functions
$\phi _{1},...,\phi_{k}\in C_{c,\mathrm{rad}}^{\infty }(B_{r_{2}}\setminus \overline{B}_{r_{1}})
$ such that $0\leq \phi _{i}\leq t_{0}$ for every $i=1,...,k$ and set 
\[
X_{k}:=\linspan\left\{ \phi _{1},...,\phi _{k}\right\} \text{\quad
and\quad }\left\| \lambda _{1}\phi _{1}+...+\lambda _{k}\phi _{k}\right\|
_{X_{k}}:=\max_{1\leq i\leq k}\left| \lambda _{i}\right| .
\]
This defines a subspace of $H_{V,\mathrm{r}}^{1}$ by assumption $\left(
h_{0}\right) $ and all norms are equivalent on $X_{k}$, so that there exist $%
m_{k},l_{k}>0$ such that for all $u\in X_{k}$ one has 
\begin{equation}
\left\| u\right\| _{X_{k}}\leq m_{k}\left\| u\right\| \quad \text{and}\quad
\left\| u\right\| _{L_{K}^{\theta }(\mathbb{R}^{N})}^{\theta }\geq l_{k}\left\|
u\right\| ^{\theta }.  \label{equivalent}
\end{equation}
Fix $\rho _{k}>0$ small enough that $km_{k}\rho _{k}<1$ and $\rho
_{k}^{2}/2-ml_{k}\rho _{k}^{\theta }<0$ (which is possible since $\theta <2$%
) and take any $u=\lambda _{1}\phi _{1}+...+\lambda _{k}\phi _{k}\in X_{k}$
such that $\left\| u\right\| =\rho _{k}$. Then by (\ref{equivalent}) we have 
\[
\left| \lambda _{i}\right| \leq \left\| u\right\| _{X_{k}}\leq m_{k}\rho
_{k}<\frac{1}{k}\qquad \text{for every }i=1,...,k
\]
and therefore 
\[
\left| u\left( x\right) \right| \leq \sum_{i=1}^{k}\left| \lambda
_{i}\right| \phi _{i}\left( x\right) \leq t_{0}\sum_{i=1}^{k}\left| \lambda
_{i}\right| <t_{0}\qquad \text{for all }x\in \mathbb{R}^{N}.
\]
By (\ref{G>}) and (\ref{equivalent}), this implies 
\[
\int_{\mathbb{R}^{N}}G\left( \left| x\right| ,u\right) dx\geq m\int_{\mathbb{R}%
^{N}}K\left( \left| x\right| \right) \left| u\right| ^{\theta }dx\geq
ml_{k}\left\| u\right\| ^{\theta }
\]
and hence we get 
\[
I\left( u\right) \leq \frac{1}{2}\left\| u\right\| ^{2}-ml_{k}\left\|
u\right\| ^{\theta }=\frac{1}{2}\rho _{k}^{2}-ml_{k}\rho _{k}^{\theta }<0.
\]
This proves (\ref{THM:wang: geom}) and the conclusion thus follows from
Theorem \ref{THM:wang}.%
\endproof

\section{Proof of the existence results for problem $\left( P_{Q}\right) $ 
\label{SEC: pfs}}

Let $\Omega \subseteq \mathbb{R}^{N}$, $N\geq 3$, be a spherically symmetric
domain (bounded or unbounded). Let $V,K,f$ be as in $\left( \mathbf{V}%
\right) ,\left( \mathbf{K}\right) ,\left( \mathbf{f}\right) $ and let $Q\in
L^{2}(\Omega _{\mathrm{r}},r^{N+1}dr)$. If $\Omega \neq \mathbb{R}^{N}$, extend
the definition of $V,K,Q$ by setting 
\[
V\left( r\right) :=+\infty \quad \text{and}\quad K\left( r\right) :=Q\left(
r\right) :=0\quad \text{for every }r\in \mathbb{R}_{+}\setminus \Omega _{%
\mathrm{r}}. 
\]
Define a Carath\'{e}odory function $g:\mathbb{R}_{+}\times \mathbb{R}\rightarrow 
\mathbb{R}$ by setting 
\[
g\left( r,t\right) :=K\left( r\right) f\left( t\right) +Q\left( r\right) . 
\]

The next lemma shows that we need only to study problem $\left( P_{Q}\right) 
$ on $\mathbb{R}^{N}$, to which the case with $\Omega \neq \mathbb{R}^{N}$
reduces. Recall the definitions (\ref{H^1_0,V :=}) and (\ref{H^1_V :=}) of
the spaces $H_{0,V}^{1}(\Omega )$ and $H_{V}^{1}(\mathbb{R}^{N})$.

\begin{lemma}
\label{PROP: space equality}If $\Omega \neq \mathbb{R}^{N}$, then, up to
restriction to $\Omega $ and null extension on $\mathbb{R}^{N}\setminus \Omega $%
, we have that $H_{V}^{1}(\mathbb{R}^{N})=H_{0,V}^{1}(\Omega )$ and any weak
solution to problem $\left( P_{Q}\right) $ on $\mathbb{R}^{N}$ is a weak
solution to $\left( P_{Q}\right) $ on $\Omega $.
\end{lemma}

Note that the result is obvious if $\Omega =\mathbb{R}^{N}$ (since $D^{1,2}(%
\mathbb{R}^{N})=D_{0}^{1,2}(\mathbb{R}^{N})$).\bigskip

\proof 
Let $u\in H_{V}^{1}(\mathbb{R}^{N})$. By the Lebesgue integration theory of
functions with real extended values (see e.g. \cite{Jones}), $\int_{\mathbb{R}%
^{N}}V\left( \left| x\right| \right) u^{2}dx<\infty $ implies $u=0$ almost
everywhere on $\mathbb{R}^{N}\setminus \Omega $ (where $V\left( \left| x\right|
\right) =+\infty $) and 
\begin{equation}
\int_{\mathbb{R}^{N}}V\left( \left| x\right| \right) u^{2}dx=\int_{\Omega
}V\left( \left| x\right| \right) u^{2}dx.  \label{0*infty}
\end{equation}
Hence $u\in D^{1,2}(\mathbb{R}^{N})$ implies 
%
%
$u\in D_{0}^{1,2}(\Omega)$ and therefore $u\in H_{0,V}^{1}(\Omega)$. 
Conversely, if $u\in H_{0,V}^{1}(\Omega)$, then $u\in D^{1,2}(\mathbb{R}%
^{N})$ and one has $\int_{\Omega }V\left( \left| x\right| \right)
u^{2}dx=\int_{\mathbb{R}^{N}}V\left( \left| x\right| \right) u^{2}\chi _{\Omega
}dx$ with $u^{2}\chi _{\Omega }=u^{2}$ almost everywhere (recall from the
statement of the lemma that we extend $u=0$ on $\mathbb{R}^{N}\setminus \Omega $%
), so that (\ref{0*infty}) holds and therefore $u\in H_{V}^{1}(\mathbb{R}^{N})$%
. This proves that $H_{V}^{1}(\mathbb{R}^{N})=H_{0,V}^{1}(\Omega )$ and the
last part of the lemma readily follows, since all the integrals involved in
the definition of weak solutions to $\left( P_{Q}\right) $ on $\mathbb{R}^{N}$
are computed on functions that vanish almost everywhere on $\mathbb{R}^{N}\setminus \Omega $.%
\endproof
\bigskip

The proof of our existence results for problem $\left( P_{Q}\right) $ relies
on the application of the general results of Section \ref{SEC: gen-results},
whose assumptions $\left( h_{0}\right) $-$\left( h_{3}\right) $ are
satisfied. Indeed, since $f\left( 0\right) =0$, one has $g\left( \cdot
,0\right) =Q$ and therefore $g$ trivially satifies assumption $\left(
h_{1}\right) $. Moreover, the potentials $V,K$ satisfy assumptions $\left(
h_{0}\right) ,\left( h_{2}\right) $ thanks to $\left( \mathbf{V}\right)
,\left( \mathbf{K}\right) $. Finally, assumption $\left( h_{3}\right) $
holds because $Q\in L^{2}(\mathbb{R}_{+},r^{N+1}dr)$ means $Q\left( \left|
\cdot \right| \right) \in L^{2}(\mathbb{R}^{N},\left| x\right| ^{2}dx)$ and
thus implies 
\[
\left| \int_{\mathbb{R}^{N}}Q\left( \left| x\right| \right) u\,dx\right| \leq
\left( \int_{\mathbb{R}^{N}}\left| x\right| ^{2}Q\left( \left| x\right| \right)
^{2}dx\right) ^{\frac{1}{2}}\left( \int_{\mathbb{R}^{N}}\frac{u^{2}}{\left|
x\right| ^{2}}dx\right) ^{\frac{1}{2}}\leq \frac{2}{N-2}\left( \int_{\mathbb{R}%
^{N}}\left| x\right| ^{2}Q\left( \left| x\right| \right) ^{2}dx\right) ^{%
\frac{1}{2}}\left\| u\right\| 
\]
for all $u\in H_{V}^{1}\hookrightarrow D^{1,2}(\mathbb{R}^{N})$, by H\"{o}lder
and Hardy inequalities.

In order to apply the general results, we will need a compactness lemma from 
\cite{BGR-p1}, which gives sufficient conditions in order that $\left( 
\mathcal{S}_{q_{1},q_{2}}^{\prime \prime }\right) $ and $\left( \mathcal{R}%
_{q_{1},q_{2}}\right) $ hold. For stating this lemma, we introduce some
functions, whose graphs are partially sketched in Figures 1-8 with a view to
easing the application of the lemma itself. For $\alpha \in \mathbb{R}$, $\beta
\in \left[ 0,1\right] $ and $\gamma \geq 2$, we define 
\[
\underline{\alpha }_{0}\left( \beta ,\gamma \right) :=\left\{ 
\begin{array}{lll}
\max \left\{ \alpha _{2}\left( \beta \right) ,\alpha _{3}\left( \beta
,\gamma \right) \right\} & \smallskip & \text{if }2\leq \gamma <N \\ 
\alpha _{1}\left( \beta ,\gamma \right) & \smallskip & \text{if }N\leq
\gamma \leq 2N-2 \\ 
-\infty &  & \text{if }\gamma >2N-2,
\end{array}
\right. 
\]
\[
\underline{q}_{0}\left( \alpha ,\beta ,\gamma \right) :=\left\{ 
\begin{array}{lll}
\max \left\{ 1,2\beta \right\} & \smallskip & \text{if }2\leq \gamma \leq N
\\ 
\max \left\{ 1,2\beta ,q_{*}\left( \alpha ,\beta ,\gamma \right) \right\} & 
\smallskip & \text{if }N<\gamma \leq 2N-2 \\ 
\max \left\{ 1,2\beta ,q_{*}\left( \alpha ,\beta ,\gamma \right)
,q_{**}\left( \alpha ,\beta ,\gamma \right) \right\} &  & \text{if }\gamma
>2N-2
\end{array}
\right. 
\]
and 
\[
\overline{q}_{0}\left( \alpha ,\beta ,\gamma \right) :=\left\{ 
\begin{array}{lll}
\min \left\{ q_{*}\left( \alpha ,\beta ,\gamma \right) ,q_{**}\left( \alpha
,\beta ,\gamma \right) \right\} & \smallskip & \text{if }2\leq \gamma <N \\ 
q_{**}\left( \alpha ,\beta ,\gamma \right) & \smallskip & \text{if }N\leq
\gamma <2N-2 \\ 
+\infty &  & \text{if }\gamma \geq 2N-2
\end{array}
\right. 
\]
(recall the definitions (\ref{alphas :=}), (\ref{q* :=}) and (\ref{q** :=})
of $\alpha _{1},\alpha _{2},\alpha _{3}$ and $q_{*},q_{**}$). For $\alpha
\in \mathbb{R}$, $\beta \in \left[ 0,1\right] $ and $\gamma \leq 2$, we define
the function 
\[
\underline{q}_{\infty }\left( \alpha ,\beta ,\gamma \right) :=\max \left\{
1,2\beta ,q_{*}\left( \alpha ,\beta ,\gamma \right) ,q_{**}\left( \alpha
,\beta ,\gamma \right) \right\} . 
\]

\begin{lemma}[{\cite[Theorems 2, 3, 4, 5]{BGR-p1}}]
\label{LEM:compactness} 
\emph{(i)} Suppose that $\Omega $ contains a neighbourhood of the origin and
assume that $\left( \mathbf{VK}_{0}\right) $ holds with $\alpha _{0}>%
\underline{\alpha }_{0}$, where $\underline{\alpha }_{0}=\underline{\alpha }%
_{0}\left( \beta _{0},2\right) $. Then 
\[
\lim_{R\rightarrow 0^{+}}\mathcal{R}_{0}\left( q_{1},R\right) =0\quad \text{%
for every }\underline{q}_{0}<q_{1}<\overline{q}_{0},
\]
where $\underline{q}_{0}=\underline{q}_{0}\left( \alpha _{0},\beta
_{0},2\right) $ and $\overline{q}_{0}=\overline{q}_{0}\left( \alpha
_{0},\beta _{0},2\right) $. If $V$ also satisfies $\left( \mathbf{V}%
_{0}\right) $, then the same result holds true with $\underline{\alpha }_{0}=%
\underline{\alpha }_{0}\left( \beta _{0},\gamma _{0}\right) $, $\underline{q}%
_{0}=\underline{q}_{0}\left( \alpha _{0},\beta _{0},\gamma _{0}\right) $ and 
$\overline{q}_{0}=\overline{q}_{0}\left( \alpha _{0},\beta _{0},\gamma
_{0}\right) $.\smallskip 

\noindent \emph{(ii)} Suppose that $\Omega $ contains a neighbourhood of
infinity and assume $\left( \mathbf{VK}_{\infty }\right) $. Then 
\[
\lim_{R\rightarrow +\infty }\mathcal{R}_{\infty }\left( q_{2},R\right)
=0\quad \text{for every }q_{2}>\underline{q}_{\infty },
\]
where $\underline{q}_{\infty }=\underline{q}_{\infty }\left( \alpha _{\infty
},\beta _{\infty },2\right) $. If $V$ also satisfies $\left( \mathbf{V}%
_{\infty }\right) $, then one can take $\underline{q}_{\infty }=\underline{q}%
_{\infty }\left( \alpha _{\infty },\beta _{\infty },\gamma _{\infty }\right) 
$.\bigskip 
\end{lemma}

\noindent 
\begin{tabular}[t]{l}
\begin{tabular}{l}
\includegraphics[height=1.5264in]{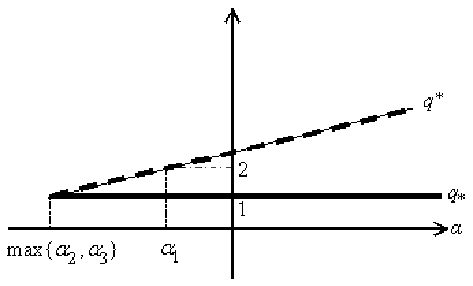}
\\ 
\textbf{Fig.1.} $\underline{q}_{0}(\cdot ,\beta ,\gamma )$ (solid) and 
$\overline{q}_{0}(\cdot ,\beta ,\gamma )$ (dashed) \\ 
for $\alpha >\underline{\alpha }_{0}$ and $\beta \in \left[ 0,1\right] $, $\gamma =2$ fixed
\end{tabular}
\end{tabular}
\begin{tabular}[t]{l}
\begin{tabular}{l}
\includegraphics[height=1.6284in]{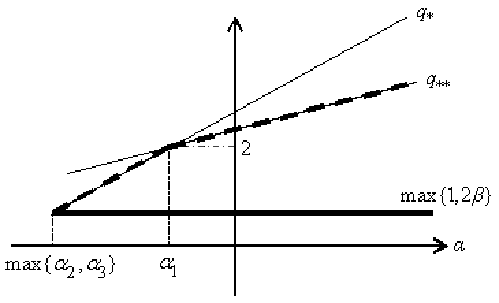}
\\ 
\textbf{Fig.2.} $\underline{q}_{0}(\cdot ,\beta ,\gamma )$ (solid) and 
$\overline{q}_{0}(\cdot ,\beta ,\gamma )$ (dashed) \\ 
for $\alpha >\underline{\alpha }_{0}$ and $\beta \in \left[ 0,1\right] $, $2\leq \gamma <N$ fixed
\end{tabular}
\end{tabular}
\bigskip

\noindent 
\begin{tabular}[t]{l}
\begin{tabular}{l}
\includegraphics[height=1.5558in]{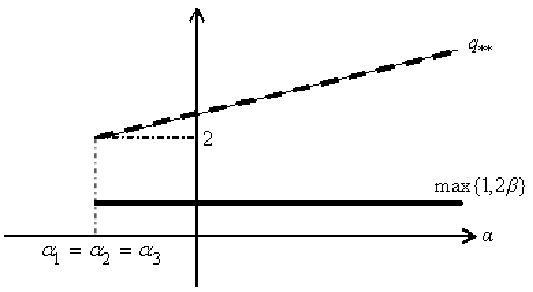}
 \\ 
\textbf{Fig.3.} $\underline{q}_{0}(\cdot ,\beta ,\gamma )$ (solid) and 
$\overline{q}_{0}(\cdot ,\beta ,\gamma )$ (dashed) \\ 
for $\alpha >\underline{\alpha }_{0}$ and $\beta \in \left[ 0,1\right] $, $\gamma =N$ fixed
\end{tabular}
\end{tabular}
\begin{tabular}[t]{l}
\begin{tabular}{l}
\includegraphics[height=1.5558in]{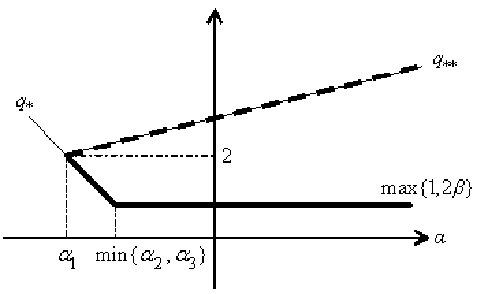}
 \\ 
\textbf{Fig.4.} $\underline{q}_{0}(\cdot ,\beta ,\gamma )$ (solid) and 
$\overline{q}_{0}(\cdot ,\beta ,\gamma )$ (dashed) \\ 
for $\alpha >\underline{\alpha }_{0}$ and $\beta \in \left[ 0,1\right] $, $N<\gamma <2N-2$ fixed
\end{tabular}
\end{tabular}
\bigskip

\noindent 
\begin{tabular}[t]{l}
\begin{tabular}{l}
\includegraphics[height=1.5558in]{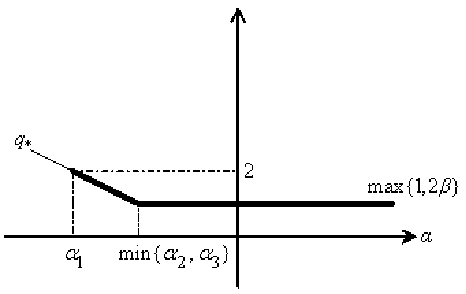}
\\ 
\textbf{Fig.5.} $\underline{q}_{0}(\cdot ,\beta ,\gamma )$ (solid) for 
$\alpha >\underline{\alpha }_{0}$ \\ 
and $\beta \in \left[ 0,1\right] $, $\gamma =2N-2$ fixed
\end{tabular}
\end{tabular}
\,\qquad 
\begin{tabular}[t]{l}
\begin{tabular}{l}
\includegraphics[height=1.5705in]{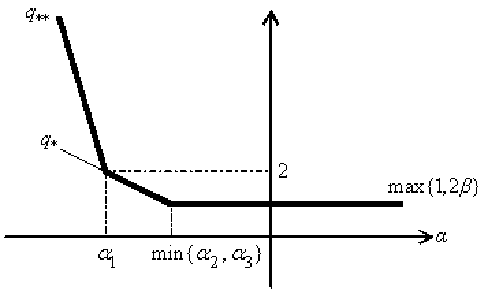}
 \\ 
\textbf{Fig.6.} $\underline{q}_{0}(\cdot ,\beta ,\gamma )$ (solid) for 
$\beta \in \left[ 0,1\right] $ \\ and $\gamma >2N-2$ fixed
\end{tabular}
\end{tabular}
\bigskip

\noindent 
\begin{tabular}[t]{l}
\begin{tabular}{l}
\includegraphics[height=1.5558in]{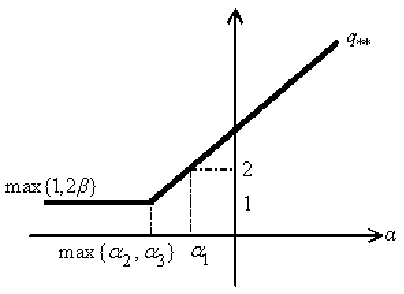}
 \\ 
\textbf{Fig.7.} $\underline{q}_{\infty }(\cdot ,\beta ,\gamma )$ (solid) for 
$\beta \in \left[ 0,1\right] $ \\ and $\gamma =2$ fixed
\end{tabular}
\end{tabular}
\,\qquad \quad \,\quad 
\begin{tabular}[t]{l}
\begin{tabular}{l}
\includegraphics[height=1.5921in]{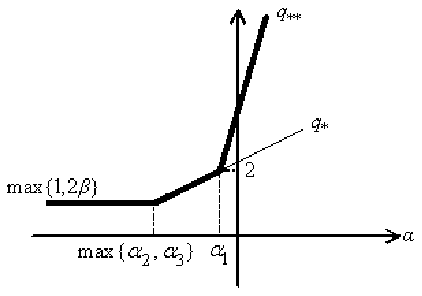}
 \\ 
\textbf{Fig.8.} $\underline{q}_{\infty }(\cdot ,\beta ,\gamma )$ (solid) for 
$\beta \in \left[ 0,1\right] $ \\ and $\gamma <2$ fixed
\end{tabular}
\end{tabular}
\bigskip

\noindent \textbf{Proof of Theorem \ref{THM:RNsuper}.}\quad 
The theorem concerns problem $\left( P_{0}\right) $, so we have $g\left( \cdot ,0\right)
=Q=0$. We want to apply Theorem \ref{THM:ex}. To this aim, since $\left( 
\mathbf{F}_{1}\right) $ implies $\left( g_{1}\right) $-$\left( g_{2}\right) $
and $\left( \mathbf{F}_{2}\right) $ implies $\left( g_{3}\right) $, we need
only to check that $\exists q_{1},q_{2}>2$ such that $\left( \mathcal{S}%
_{q_{1},q_{2}}^{\prime \prime }\right) $ and $\left( f_{q_{1},q_{2}}\right) $
hold. This follows from the other assumptions of the theorem, since $\left( 
\mathbf{f}_{1}\right) $ and $f\left( 0\right) =0$ imply $\left(
f_{q_{1},q_{2}}\right) $ for $t\geq 0$ (which is actually enough by Remark 
\ref{RMK:thm:ex}.\ref{RMK:thm:ex-2}), and inequalities (\ref{THM:RNsuper:
ineqs}) imply $q_{1},q_{2}>2$ and 
\[
\lim_{R\rightarrow 0^{+}}\mathcal{R}_{0}\left( q_{1},R\right)
=\lim_{R\rightarrow +\infty }\mathcal{R}_{\infty }\left( q_{2},R\right) =0 
\]
by Lemma \ref{LEM:compactness} (cf. Figs. 1-8). Note that also Proposition 
\ref{PROP:symm-crit} applies, since $\left( \mathcal{R}_{q_{1},q_{2}}\right) 
$ holds. Hence the proof is complete, as the result follows from Theorem \ref
{THM:ex} and Proposition \ref{PROP:symm-crit}.%
\endproof
\bigskip

\noindent \textbf{Proof of Theorems \ref{THM:Bsuper} and \ref{THM:BCsuper}.}\quad 
The proof runs exactly as the proof of Theorem \ref{THM:RNsuper} and
then ends by applying Lemma \ref{PROP: space equality}. The only difference
is that, under the assumptions of Theorem \ref{THM:Bsuper}, we use
inequalities (\ref{THM:Bsuper: ineqs}) to deduce $q>2$ and 
\[
\lim_{R\rightarrow 0^{+}}\mathcal{R}_{0}\left( q,R\right) =0 
\]
by part (i) of Lemma \ref{LEM:compactness} (cf. Figs. 1-6), and then we
observe that $\mathcal{R}_{\infty }\left( q,R\right) =0$ for every $R>0$
large enough, since the integrand $K\left( \left| \cdot \right| \right)
\left| u\right| ^{q-1}\left| h\right| $ vanishes almost everywhere outside $%
\Omega $. Similarly, under the assumptions of Theorem \ref{THM:BCsuper}, we
deduce 
\[
\lim_{R\rightarrow +\infty }\mathcal{R}_{\infty }\left( q_{2},R\right) =0 
\]
by part (ii) of Lemma \ref{LEM:compactness} (cf. Figs. 7-8), and then we
observe that $\mathcal{R}_{0}\left( q,R\right) =0$ for every $R>0$ small
enough. Hence, in both cases, conditions $\left( \mathcal{S}%
_{q_{1},q_{2}}^{\prime \prime }\right) $, $\left( \mathcal{R}%
_{q_{1},q_{2}}\right) $ and $\left( f_{q_{1},q_{2}}\right) $ hold with $%
q_{1}=q_{2}=q$ and we can apply Theorem \ref{THM:ex} and Proposition \ref
{PROP:symm-crit}.%
\endproof
\bigskip

\noindent \textbf{Proof of Theorems \ref{THM:RNsub1} and \ref{THM:RNsub2}.}\quad 
The proof is analogous to the one of Theorem \ref{THM:RNsuper} and
reduces to the application of Corollary \ref{COR:ex-sub}. We have $g\left(
\cdot ,0\right) =Q\geq 0$ by assumption. This implies $G\left( r,t\right)
=K\left( r\right) F\left( t\right) +Q\left( r\right) t\geq K\left( r\right)
F\left( t\right) $ for almost every $r>0$ and all $t\geq 0$, so that $\left(
g_{6}\right) $ or $\left( g_{7}\right) $ holds true, respectively according
as $f$ satisfies $\left( \mathbf{F}_{3}\right) $ or $Q$ does not vanish
almost everywhere in $\left( r_{1},r_{2}\right) $. Thus, in order to apply
Corollary \ref{COR:ex-sub}, we need only to check that the other assumptions
of the theorems yield the existence of $q_{1},q_{2}\in \left( 1,2\right) $
such that $\left( \mathcal{S}_{q_{1},q_{2}}^{\prime \prime }\right) $ and $%
\left( f_{q_{1},q_{2}}\right) $ hold. On the one hand, inequalities (\ref
{THM:RNsub1: ineqs}) and (\ref{THM:RNsub2: ineqs}) imply $q,q_{1},q_{2}\in
\left( 1,2\right) $ and 
\[
\lim_{R\rightarrow 0^{+}}\mathcal{R}_{0}\left( q,R\right)
=\lim_{R\rightarrow +\infty }\mathcal{R}_{\infty }\left( q,R\right)
=\lim_{R\rightarrow 0^{+}}\mathcal{R}_{0}\left( q_{1},R\right)
=\lim_{R\rightarrow +\infty }\mathcal{R}_{\infty }\left( q_{2},R\right) =0 
\]
by Lemma \ref{LEM:compactness} (cf. Figs. 1-8). This also allows us to apply
Proposition \ref{PROP:symm-crit}, since $\left( \mathcal{R}%
_{q_{1},q_{2}}\right) $ holds. On the other hand, $\left( \mathbf{f}%
_{1}\right) $ and $\left( \mathbf{f}_{2}\right) $ respectively imply $\left(
f_{q_{1},q_{2}}\right) $ and $\left( f_{q,q}\right) $ for $t\geq 0$, which
is actually enough by Remark \ref{RMK:ex-sub}.\ref{RMK:ex-sub-2}. Hence the
result follows from Corollary \ref{COR:ex-sub} and Proposition \ref
{PROP:symm-crit}, so that the proof is complete.%
\endproof
\bigskip

The proof of Theorem \ref{THM:multiple} is very similar to the ones of
Theorems \ref{THM:RNsuper}, \ref{THM:RNsub1} and \ref{THM:RNsub2}, so we
leave it to the interested reader (apply Theorems \ref{THM:mult} and 
\ref{THM:mult-sub} instead of Theorem \ref{THM:ex} and Corollary \ref{COR:ex-sub}).

\section{Appendix}

This Appendix is devoted to the derivation of some radial estimates, which
has been announced and used in \cite[Lemmas 3 and 4]{BGR-p1} to prove the
compactness results yielding Lemma \ref{LEM:compactness} above. Such
estimates has been also proved in \cite[Lemmas 4 and 5]{Su-Wang-Will-p}, but
under slightly different assumptions from the ones used here (and in \cite
{BGR-p1}).

Assume $N\geq 3$ and denote by $\sigma _{N}$ the $(N-1)$-dimensional measure
of the unit sphere of $\mathbb{R}^{N}$.

\begin{lemma}
\label{Lem: W(J)}Let $u$ be any radial map of $D^{1,2}(\mathbb{R}^{N})$ and let 
$\tilde{u}:\left( 0,+\infty \right) \rightarrow \mathbb{R}$ be such that $%
u\left( x\right) =\tilde{u}\left( \left| x\right| \right) $ for almost every 
$x\in \mathbb{R}^{N}$. Then $\tilde{u}\in W^{1,1}\left( \mathcal{I}\right) $
for any bounded open interval $\mathcal{I}\subset \left( 0,+\infty \right) $
such that $\inf \mathcal{I}>0$.
\end{lemma}

\proof 
Denote $\lambda :=\inf \mathcal{I}>0$ and set $\mathcal{I}_{N}:=\left\{ x\in 
\mathbb{R}^{N}:\left| x\right| \in \mathcal{I}\right\} $. Since $u\in L_{%
\mathrm{loc}}^{1}(\mathbb{R}^{N})$, we readily have 
\[
\int_{\mathcal{I}}\left| \tilde{u}\left( r\right) \right| dr\leq \frac{1}{%
\lambda ^{N-1}}\int_{\mathcal{I}}r^{N-1}\left| \tilde{u}\left( r\right)
\right| dr=\frac{1}{\lambda ^{N-1}\sigma _{N}}\int_{\mathcal{I}_{N}}\left|
u\left( x\right) \right| dx<\infty \,. 
\]
Now we exploit the density of $C_{\mathrm{c}}^{\infty }(\mathbb{R}^{N})$ radial
mappings in the space of $D^{1,2}(\mathbb{R}^{N})$ radial mappings (which
follows from standard convolution and regularization arguments) in order to
infer that, as $u$ is radial, $\nabla u\left( x\right) \cdot x$ only depends
on $\left| x\right| $. Thus there exists $\phi :\left( 0,+\infty \right)
\rightarrow \mathbb{R}$ such that $\nabla u\left( x\right) \cdot x=\phi \left(
\left| x\right| \right) $ for almost every $x\in \mathbb{R}^{N}$ and one has 
\[
\int_{\mathcal{I}}\frac{\left| \phi \left( r\right) \right| }{r}dr=\int_{%
\mathcal{I}}r^{N-1}\frac{\left| \phi \left( r\right) \right| }{r^{N}}dr=%
\frac{1}{\sigma _{N}}\int_{\mathcal{I}_{N}}\frac{\left| \nabla u\left(
x\right) \cdot x\right| }{\left| x\right| ^{N}}dx\leq \frac{1}{\lambda
^{N-1}\sigma _{N}}\int_{\mathcal{I}_{N}}\left| \nabla u\left( x\right)
\right| dx<\infty \,. 
\]
Then, letting $\varphi \in C_{c}^{\infty }(\mathcal{I})$ and setting $\psi
\left( x\right) :=\varphi \left( \left| x\right| \right) $, we get 
\begin{eqnarray*}
\int_{\mathcal{I}}\tilde{u}\left( r\right) \varphi ^{\prime }\left( r\right)
dr &=&\frac{1}{\sigma _{N}}\int_{\mathcal{I}_{N}}\frac{\tilde{u}\left(
\left| x\right| \right) }{\left| x\right| ^{N-1}}\varphi ^{\prime }\left(
\left| x\right| \right) dx=\frac{1}{\sigma _{N}}\int_{\mathcal{I}_{N}}\frac{%
u\left( x\right) }{\left| x\right| ^{N-1}}\sum_{i=1}^{N}\frac{\partial \psi 
}{\partial x_{i}}\left( x\right) \frac{x_{i}}{\left| x\right| }dx \\
&=&-\frac{1}{\sigma _{N}}\sum_{i=1}^{N}\int_{\mathcal{I}_{N}}\frac{\partial u%
}{\partial x_{i}}\left( x\right) \frac{x_{i}}{\left| x\right| ^{N}}\psi
\left( x\right) dx=-\frac{1}{\sigma _{N}}\int_{\mathcal{I}_{N}}\frac{\nabla
u\left( x\right) \cdot x}{\left| x\right| ^{N}}\psi \left( x\right) dx \\
&=&-\frac{1}{\sigma _{N}}\int_{\mathcal{I}_{N}}\frac{\phi \left( \left|
x\right| \right) }{\left| x\right| ^{N}}\varphi \left( \left| x\right|
\right) dx=-\int_{\mathcal{I}}\frac{\phi \left( r\right) }{r}\varphi \left(
r\right) dr
\end{eqnarray*}
because $\psi \in C_{c}^{\infty }(\mathcal{I}_{N})$.%
\endproof
\bigskip

Let $V:\mathbb{R}_{+}\rightarrow \left[ 0,+\infty \right] $ be a measurable
function satisfying $\left( h_{0}\right) $ and define $H_{V,\mathrm{r}}^{1}$
as in (\ref{H^1_Vr :=}).

\begin{proposition}
\label{A: stima1}Assume that there exists $R_{2}>0$ such that $V\left(
r\right) <+\infty $ for almost every $r>R_{2}$ and 
\[
\lambda _{\infty }:=\essinf_{r>R_{2}}r^{\gamma _{\infty }}V\left(
r\right) >0\quad \text{for some }\gamma _{\infty }\leq 2.
\]
Then every $u\in H_{V,\mathrm{r}}^{1}$ satisfies 
\[
\left| u\left( x\right) \right| \leq c_{\infty }\lambda _{\infty }^{-\frac{1%
}{4}}\left\| u\right\| \left| x\right| ^{-\frac{2(N-1)-\gamma _{\infty }}{4}%
}\quad \text{almost everywhere in }B_{R_{2}}^{c},
\]
where $c_{\infty }=\sqrt{2/\sigma _{N}}$.
\end{proposition}

\proof 
Let $u\in H_{V,\mathrm{r}}^{1}$ and let $\tilde{u}:\left( 0,+\infty \right)
\rightarrow \mathbb{R}$ be continuous and such that $u\left( x\right) =\tilde{u}%
\left( \left| x\right| \right) $ for almost every $x\in \mathbb{R}^{N}$. Set 
\[
v\left( r\right) :=r^{N-1-\gamma _{\infty }/2}\tilde{u}\left( r\right)
^{2}\quad \text{for all }r>0.
\]
If $\lambda :=\liminf\limits_{r\rightarrow +\infty }v\left( r\right) >0$,
then for every $r$ large enough one has 
\[
r^{N-1-\gamma _{\infty }}\tilde{u}\left( r\right) ^{2}\geq \frac{\lambda }{%
2r^{\gamma _{\infty }/2}},
\]
whence, since $\gamma _{\infty }\leq 2$, one gets the contradiction 
\[
\int_{B_{R_{2}}^{c}}V\left( \left| x\right| \right) u^{2}dx\geq \lambda
_{\infty }\int_{B_{R_{2}}^{c}}\frac{u^{2}}{\left| x\right| ^{\gamma _{\infty
}}}dx=\lambda _{\infty }\sigma _{N}\int_{R_{2}}^{+\infty }\frac{\tilde{u}%
\left( r\right) ^{2}}{r^{\gamma _{\infty }}}r^{N-1}dr\geq \lambda _{\infty
}\sigma _{N}\int_{R_{2}}^{+\infty }\frac{\lambda }{2r^{\gamma _{\infty }/2}}%
dr=+\infty .
\]
Therefore it must be $\lambda =0$ and thus there exists $r_{n}\rightarrow
+\infty $ such that $v\left( r_{n}\right) \rightarrow 0$. By Lemma \ref{Lem:
W(J)}, we have $v\in W^{1,1}\left( \left( r,r_{n}\right) \right) $ for every 
$R_{2}<r<r_{n}<+\infty $ and hence 
\[
v\left( r_{n}\right) -v\left( r\right) =\int_{r}^{r_{n}}v^{\prime }\left(
s\right) ds\,.
\]
Moreover, for almost every $s\in \left( r,r_{n}\right) $ one has 
\begin{eqnarray*}
v^{\prime }\left( s\right)  &=&\left( N-1-\frac{\gamma _{\infty }}{2}\right)
s^{N-2-\frac{\gamma _{\infty }}{2}}\tilde{u}\left( s\right) ^{2}+2s^{N-1-%
\frac{\gamma _{\infty }}{2}}\tilde{u}\left( s\right) \tilde{u}^{\prime
}\left( s\right) \geq 2s^{N-1-\frac{\gamma _{\infty }}{2}}\tilde{u}\left(
s\right) \tilde{u}^{\prime }\left( s\right)  \\
&\geq &-2s^{N-1-\frac{\gamma _{\infty }}{2}}\left| \tilde{u}\left( s\right)
\right| \left| \tilde{u}^{\prime }\left( s\right) \right| 
\end{eqnarray*}
(note that $N-1-\frac{\gamma _{\infty }}{2}\geq N-2>0$). Hence 
\[
v\left( r_{n}\right) -v\left( r\right) =\int_{r}^{r_{n}}v^{\prime }\left(
s\right) ds\geq -2\int_{r}^{r_{n}}s^{N-1-\frac{\gamma _{\infty }}{2}}\left| 
\tilde{u}\left( s\right) \right| \left| \tilde{u}^{\prime }\left( s\right)
\right| ds
\]
and therefore 
\begin{eqnarray*}
v\left( r\right) -v\left( r_{n}\right)  &\leq &2\int_{r}^{r_{n}}s^{\frac{N-1%
}{2}}\frac{\left| \tilde{u}\left( s\right) \right| }{s^{\gamma _{\infty }/2}}%
s^{\frac{N-1}{2}}\left| \tilde{u}^{\prime }\left( s\right) \right| ds\leq
2\left( \int_{r}^{r_{n}}s^{N-1}\frac{\tilde{u}\left( s\right) ^{2}}{%
s^{\gamma _{\infty }}}ds\right) ^{\frac{1}{2}}\left( \int_{r}^{r_{n}}s^{N-1}%
\tilde{u}^{\prime }\left( s\right) ^{2}ds\right) ^{\frac{1}{2}} \\
&\leq &2\left( \int_{R_{2}}^{+\infty }s^{N-1}\frac{\tilde{u}\left( s\right)
^{2}}{s^{\gamma _{\infty }}}ds\right) ^{\frac{1}{2}}\left( \int_{0}^{+\infty
}s^{N-1}\tilde{u}^{\prime }\left( s\right) ^{2}ds\right) ^{\frac{1}{2}} \\
&\leq &2\left( \frac{1}{\lambda _{\infty }}\int_{R_{2}}^{+\infty }V\left(
s\right) \tilde{u}\left( s\right) ^{2}s^{N-1}ds\right) ^{\frac{1}{2}}\left(
\int_{0}^{+\infty }\tilde{u}^{\prime }\left( s\right) ^{2}s^{N-1}ds\right) ^{%
\frac{1}{2}} \\
&\leq &\frac{2}{\sigma _{N}}\left( \frac{1}{\lambda _{\infty }}\int_{\mathbb{R}%
^{N}}V\left( \left| x\right| \right) u^{2}dx\right) ^{\frac{1}{2}}\left(
\int_{\mathbb{R}^{N}}\left| \nabla u\right| ^{2}dx\right) ^{\frac{1}{2}}.
\end{eqnarray*}
This implies $v\left( r\right) \leq 2\lambda _{\infty }^{-1/2}\left\|
u\right\| ^{2}/\sigma _{N}$, whence the result readily ensues.%
\endproof
\bigskip 

\begin{proposition}
\label{A: stima2}Assume that there exists $R>0$ such that $V\left( r\right)
<+\infty $ almost everywhere on $\left( 0,R\right) $ and 
\[
\lambda _{0}:=\essinf_{r\in \left( 0,R\right) }r^{\gamma
_{0}}V\left( r\right) >0\quad \text{for some }\gamma _{0}\geq 2.
\]
Then every $u\in H_{V,\mathrm{r}}^{1}$ satisfies 
\[
\left| u\left( x\right) \right| \leq c_{0}\left( \frac{1}{\sqrt{\lambda _{0}}%
}+\frac{R^{\frac{\gamma _{0}-2}{2}}}{\lambda _{0}}\right) ^{\frac{1}{2}%
}\left\| u\right\| \left| x\right| ^{-\frac{2N-2-\gamma _{0}}{4}}\quad \text{%
almost everywhere in }B_{R},
\]
where $c_{0}=\sqrt{\max \left\{ 2,N-2\right\} /\sigma _{N}}$.
\end{proposition}

\proof 
Let $u\in H_{V,\mathrm{r}}^{1}$ and $\tilde{u}:\left( 0,+\infty \right)
\rightarrow \mathbb{R}$ continuous and such that $u\left( x\right) =\tilde{u}%
\left( \left| x\right| \right) $ for almost every $x\in \mathbb{R}^{N}$. Set 
\[
v\left( r\right) :=r^{N-1-\gamma _{0}/2}\tilde{u}\left( r\right) ^{2}\quad 
\text{for all }r>0.
\]
If $\lambda :=\liminf\limits_{r\rightarrow 0^{+}}v\left( r\right) >0$, then
for every $r$ small enough one gets 
\[
r^{N-1-\gamma _{0}}\tilde{u}\left( r\right) ^{2}\geq \frac{\lambda }{%
2r^{\gamma _{0}/2}}
\]
and thus, since $\gamma _{0}\geq 2$, one deduces the contradiction 
\[
\int_{B_{R}}V\left( \left| x\right| \right) u^{2}dx\geq \lambda
_{0}\int_{B_{R}}\frac{u^{2}}{\left| x\right| ^{\gamma _{0}}}dx=\lambda
_{0}\sigma _{N}\int_{0}^{R}\frac{\tilde{u}\left( r\right) ^{2}}{r^{\gamma
_{0}}}r^{N-1}dr\geq \lambda _{0}\sigma _{N}\int_{0}^{R}\frac{\lambda }{%
2r^{\gamma _{0}/2}}dr=+\infty .
\]
So we have $\lambda =0$ and thus there exists $r_{n}\rightarrow 0^{+}$ such
that $v\left( r_{n}\right) \rightarrow 0$. By Lemma \ref{Lem: W(J)}, we have 
$v\in W^{1,1}\left( \left( r_{n},r\right) \right) $ for every $0<r_{n}<r<R$
and hence 
\[
v\left( r\right) -v\left( r_{n}\right) =\int_{r_{n}}^{r}v^{\prime }\left(
s\right) ds\,.
\]
Moreover, for almost every $s\in \left( r_{n},r\right) $ one has 
\[
v^{\prime }\left( s\right) =\left( N-1-\frac{\gamma _{0}}{2}\right) s^{N-2-%
\frac{\gamma _{0}}{2}}\tilde{u}\left( s\right) ^{2}+2s^{N-1-\frac{\gamma _{0}%
}{2}}\tilde{u}\left( s\right) \tilde{u}^{\prime }\left( s\right) .
\]
Notice that 
\begin{eqnarray*}
\int_{r_{n}}^{r}s^{N-2-\frac{\gamma _{0}}{2}}\tilde{u}\left( s\right) ^{2}ds
&=&\int_{r_{n}}^{r}s^{\frac{\gamma _{0}-2}{2}}\frac{\tilde{u}\left( s\right)
^{2}}{s^{\gamma _{0}}}s^{N-1}ds\leq r^{\frac{\gamma _{0}-2}{2}}\int_{0}^{R}%
\frac{\tilde{u}\left( s\right) ^{2}}{s^{\gamma _{0}}}s^{N-1}ds \\
&\leq &r^{\frac{\gamma _{0}-2}{2}}\frac{1}{\lambda _{0}}\int_{0}^{R}V\left(
s\right) \tilde{u}\left( s\right) ^{2}s^{N-1}ds\leq R^{\frac{\gamma _{0}-2}{2%
}}\frac{1}{\lambda _{0}}\frac{1}{\sigma _{N}}\left\| u\right\| ^{2}
\end{eqnarray*}
and 
\begin{eqnarray*}
\int_{r_{n}}^{r}s^{N-1-\frac{\gamma _{0}}{2}}\tilde{u}\left( s\right) 
\tilde{u}^{\prime }\left( s\right) ds &\leq &\int_{r_{n}}^{r}s^{N-1-\frac{%
\gamma _{0}}{2}}\left| \tilde{u}\left( s\right) \right| \left| \tilde{u}%
^{\prime }\left( s\right) \right| ds \\
&\leq &\left( \int_{r_{n}}^{r}s^{N-1}\frac{\tilde{u}\left( s\right) ^{2}}{%
s^{\gamma _{0}}}ds\right) ^{\frac{1}{2}}\left( \int_{r_{n}}^{r}s^{N-1}%
\tilde{u}^{\prime }\left( s\right) ^{2}ds\right) ^{\frac{1}{2}} \\
&\leq &\left( \int_{0}^{R}s^{N-1}\frac{\tilde{u}\left( s\right) ^{2}}{%
s^{\gamma _{0}}}ds\right) ^{\frac{1}{2}}\left( \int_{0}^{+\infty }s^{N-1}%
\tilde{u}^{\prime }\left( s\right) ^{2}ds\right) ^{\frac{1}{2}} \\
&\leq &\left( \frac{1}{\lambda _{0}}\int_{0}^{R}V\left( s\right) \tilde{u}%
\left( s\right) ^{2}s^{N-1}ds\right) ^{\frac{1}{2}}\left( \int_{0}^{+\infty }%
\tilde{u}^{\prime }\left( s\right) ^{2}s^{N-1}ds\right) ^{\frac{1}{2}} \\
&\leq &\frac{1}{\sigma _{N}}\left( \frac{1}{\lambda _{0}}\int_{\mathbb{R}%
^{N}}V\left( \left| x\right| \right) u^{2}dx\right) ^{\frac{1}{2}}\left(
\int_{\mathbb{R}^{N}}\left| \nabla u\right| ^{2}dx\right) ^{\frac{1}{2}}=\frac{1%
}{\sigma _{N}}\frac{1}{\lambda _{0}^{1/2}}\left\| u\right\| ^{2}.
\end{eqnarray*}
If $2\leq \gamma _{0}<2N-2$ (i.e., $N-1-\frac{\gamma _{0}}{2}>0$), then 
\begin{eqnarray*}
v\left( r\right) -v\left( r_{n}\right)  &=&\int_{r_{n}}^{r}v^{\prime }\left(
s\right) ds\leq \left( N-1-\frac{\gamma _{0}}{2}\right)
\int_{r_{n}}^{r}s^{N-2-\frac{\gamma _{0}}{2}}\tilde{u}\left( s\right)
^{2}ds+2\int_{r_{n}}^{r}s^{N-1-\frac{\gamma _{0}}{2}}\tilde{u}\left(
s\right) \tilde{u}^{\prime }\left( s\right) ds \\
&\leq &\left( N-1-\frac{\gamma _{0}}{2}\right) R^{\frac{\gamma _{0}-2}{2}}%
\frac{1}{\lambda _{0}}\frac{1}{\sigma _{N}}\left\| u\right\| ^{2}+\frac{2}{%
\sigma _{N}}\frac{1}{\lambda _{0}^{1/2}}\left\| u\right\| ^{2} \\
&\leq &\left( N-2\right) R^{\frac{\gamma _{0}-2}{2}}\frac{1}{\lambda _{0}}%
\frac{1}{\sigma _{N}}\left\| u\right\| ^{2}+\frac{2}{\sigma _{N}}\frac{1}{%
\lambda _{0}^{1/2}}\left\| u\right\| ^{2}.
\end{eqnarray*}
If $\gamma _{0}\geq 2N-2$ (i.e., $N-1-\frac{\gamma _{0}}{2}\leq 0$), then 
\[
v^{\prime }\left( s\right) =\left( N-1-\frac{\gamma _{0}}{2}\right) s^{N-2-%
\frac{\gamma _{0}}{2}}\tilde{u}\left( s\right) ^{2}+2s^{N-1-\frac{\gamma _{0}%
}{2}}\tilde{u}\left( s\right) \tilde{u}^{\prime }\left( s\right) \leq
2s^{N-1-\frac{\gamma _{0}}{2}}\tilde{u}\left( s\right) \tilde{u}^{\prime
}\left( s\right) 
\]
and therefore 
\[
v\left( r\right) -v\left( r_{n}\right) =\int_{r_{n}}^{r}v^{\prime }\left(
s\right) ds\leq 2\int_{r_{n}}^{r}s^{N-1-\frac{\gamma _{0}}{2}}\tilde{u}%
\left( s\right) \tilde{u}^{\prime }\left( s\right) ds\leq \frac{2}{\sigma
_{N}}\frac{1}{\lambda _{0}^{1/2}}\left\| u\right\| ^{2}.
\]
So, in any case, we have 
\[
v\left( r\right) -v\left( r_{n}\right) \leq \frac{1}{\sigma _{N}}\left( 
\frac{2}{\lambda _{0}^{1/2}}+\left( N-2\right) \frac{R^{\frac{\gamma _{0}-2}{%
2}}}{\lambda _{0}}\right) \left\| u\right\| ^{2},
\]
which implies 
\[
v\left( r\right) \leq \frac{1}{\sigma _{N}}\left( \frac{2}{\lambda _{0}^{1/2}%
}+\left( N-2\right) \frac{R^{\frac{\gamma _{0}-2}{2}}}{\lambda _{0}}\right)
\left\| u\right\| ^{2}.
\]
Hence the result readily follows.%
\endproof

\begin{remark}
By Hardy inequality, the radial subspace of $D^{1,2}(\mathbb{R}^{N})$ equals $%
H_{1/\left| x\right| ^{2},\mathrm{r}}^{1}$, with equivalent norms.
Therefore, taking $\gamma _{0}=\gamma _{\infty }=2$ and $R=R_{2}$ in
Propositions \ref{A: stima1} and \ref{A: stima2}, we get $\lambda
_{0}=\lambda _{\infty }=1$ and find the renowned Radial Lemma known as Ni's
inequality: there exists $c_{N}>0$ such that every radial function in $%
D^{1,2}(\mathbb{R}^{N})$ satisfies 
$\left| u\left( x\right) \right| \leq c_{N}\left\| \nabla u\right\| _{L^{2}(%
\mathbb{R}^{N})}\left| x\right| ^{-(N-2)/2}$ almost everywhere in $\mathbb{R}^{N}$.
\end{remark}

\end{document}